\documentclass{amsart} 
\usepackage{amssymb}
\usepackage{amscd}
\usepackage{graphicx}

\theoremstyle{plain}
\newtheorem{theorem}{Theorem}[section]
\newtheorem{lemma}[theorem]{Lemma}

\newtheorem{proposition}[theorem]{Proposition}

\theoremstyle{definition}
\newtheorem{definition}[theorem]{Definition}
\newtheorem{example}[theorem]{Example}
\newtheorem{examples}[theorem]{Examples}
\newtheorem{notation}[theorem]{Notation}
\newtheorem*{standing assumption}{Standing Assumption}

\theoremstyle{remark}
\newtheorem{remark}[theorem]{Remark}
\newtheorem{remarks}[theorem]{Remarks}
\newtheorem*{ack}{Acknowledgment}

\begin{document}
\title[One-dimensional solenoids]
{Canonical symbolic dynamics for\\
one-dimensional generalized solenoids}

\author{Yi, Inhyeop}

\address{Department of Mathematics, 
University of Maryland, 
College park, MD, 20742}

\email{inhyeop@math.umd.edu}



\begin{abstract}
We define canonical subshift of finite type
cover for Williams' 1-dimensional generalized
solenoids, and use resulting invariants
to distinguish some closely related solenoids.
\end{abstract}
\maketitle

\section{Introduction}\label{S1}

R.~F.~Williams  has developed
a theory of expanding attractors
for a dynamical system (\cite{w1, w2, w3}). 
These can be modeled as shift maps of
generalized $n$-solenoids which are defined as 
inverse limits of immersions of
$n$-dimensional branched manifolds
satisfying certain axioms.

In this paper, we produce canonical 
shift of finite type (SFT) covers of Williams' 
$1$-solenoids in the following sense:
Let $\overline{X}$ be a $1$-solenoid, 
$\overline f$ a shift map on $\overline{X}$, 
and $\overline{\mathcal{O}}$ 
a union of finitely many periodic orbits of 
$\overline f$. 
We give an algorithm which  takes the input
$\{\overline{X},\overline{f},\overline{\mathcal{O}}\}$ 
and produces as output a mixing SFT 
$\Sigma_{\mathcal{O}}$ with shift map 
$\sigma_{\mathcal{O}}$ and a semiconjugacy 
$p_{\mathcal{O}}\colon
\Sigma_{\mathcal{O}} \to \overline X$
(that is, $p_{\mathcal{O}}$ is a continuous surjection
and $p_{\mathcal{O}}\circ \sigma_{\mathcal{O}}
=\overline{f} \circ p_{\mathcal{O}}$).
Then we prove that if there is a conjugacy $\phi$ 
of $1$-solenoids  
$(\overline{X},\overline{f})$ and 
$(\overline{X'},\overline{f'})$ such that 
$\phi$ sends  
$\overline{\mathcal{O}}$ to $\overline{\mathcal{O'}}$, 
then there is a unique conjugacy $\tilde \phi$
such that
$\phi \circ p_{\mathcal O}
=p_{\mathcal O'}\circ \tilde{\phi}$. 
These covers can be exploited
to give nontrivial computable 
invariants distinguishing closely 
related solenoid maps. 

To our knowledge, we are giving 
the first construction of canonical SFT covers
for a class of nonzero dimensional systems. 
The canonicalness requires
the dependence on $\mathcal O$, 
and perhaps this is why it was not noticed earlier. 
However, there have been other works (\cite{am, sn})
achieving some specialness or invariant
for SFT covers of some systems,  
and there were earlier constructions of canonical 
covers for some systems. 
Krieger (\cite{kr}) gave canonical SFT covers of 
sofic systems (which are zero dimensional),
and Fried (\cite{b2, fd}) 
more generally offered canonical coordinate (CC)  
covers of finitely presented (FP) systems.  
These covers are built from
sets of possible pasts and futures.
The Krieger-Fried covers make the step from FP to CC. 
The $1$-solenoids are already CC,
and the covers we produce are SFTs.
We raise the question,
can our relatively simple one-dimensional
construction be generalized in some inductive way
to produce canonical symbolic dynamics for 
higher dimensional generalized solenoids? 
 
Apart from the matter of canonical symbolic dynamics,
we mention renewed interest in Williams' systems and
related systems on account of connections
with ordered group invariants (\cite{bjv, sv, y1}) and
substitutions and tilings (\cite{ap,dhs}). 

We study the $1$-solenoids as
purely topological systems. 
For this we give some defining topological axioms
closely related to Williams' axioms.
A $1$-solenoid of Williams becomes
one of our $1$-solenoids by ignoring
the differentiable structure.
Conversely, every topologically defined $1$-solenoid
can be given a differentiable structure 
which makes it a $1$-solenoid in the sense of Williams. 
However, the essential aspects of the situation
are not differentiable but topological,
and to clarify this
we give the purely topological development.

The outline of the paper is as follows.
In section 2,
following Williams (\cite{w1, w2}) rather closely,
we give axioms for our systems and
prove some basic facts about them.
We also recall the construction of 
an SFT cover from a graph presentation. 
In section 3, we recall Williams' definition of 
shift equivalence and show that
every topological conjugacy of branched $1$-solenoids
is induced by a shift equivalence of
their graph presentations. 
(A `branched' solenoid is
a solenoid derived from a presentation
which need not satisfy Williams' Flattening Axiom,
so this is a slight generalization of Williams' work.) 
We also establish a key observation:
If the shift equivalence is given by graph maps
(maps sending vertices to vertices), 
then the conjugacy lifts uniquely
to a conjugacy of the SFT covers derived from 
the graph presentations. 
In section 4, given $\overline{\mathcal O}$, 
we give a graph algorithm for
a new graphical presentation
$(X_{\mathcal{O}},f_{\mathcal{O}})$
of  the solenoid system. 
One consequence of this construction is that
every $1$-solenoid with a fixed point has
an elementary presentation in the sense of Williams,
so this extends
Williams' classification result (\cite[\S7]{w2})
to all $1$-solenoids with fixed points.
In Section 5, we use the previous results to 
produce the canonical SFT covers, and use them 
to distinguish the pair of systems considered 
by Williams and Ustinov (\cite{us, w2})
by computing Bowen-Franks groups of
certain attached canonical SFT covers.   
In Appendix A, 
we show that our canonical 
SFT covers are not canonical as one-sided SFTs, 
despite the one-sided aspects of the construction.  
In Appendix B, we show that our topological 
$1$-solenoids can be given a differentiable structure
making them differentiable immersions in the sense 
of Williams.

\section{Markov maps and their SFT covers}\label{S2}

In the style of Williams (\cite{w1, w2}),
we will define several axioms which 
might be satisfied by a continuous self-map
of a directed graph.
Let $X$ be a directed graph with vertex set
$\mathcal{V}$ and edge set $\mathcal{E}$,
and $f\colon X\to X$ a continuous map.
Axioms $0$-$3$ and $5$ correspond to Williams'
Axioms $0$-$2$, $3^\circ$, and $4$ in \cite{w2}.
\begin{enumerate}
\item[Axiom 0.]({\it Indecomposability})
$(X,f)$ is indecomposable.
\item[Axiom 1.]({\it Nonwandering})
All points of $X$ are nonwandering under $f$.
\item[Axiom 2.] ({\it Flattening}) 
There is $k\ge 1$ such that  for all $x\in X$ 
there is an open neighborhood $U$ of $x$
such that $f^k(U)$ is homeomorphic to
$(-\epsilon, \epsilon )$.
\item[Axiom 3.]({\it Expansion}) 
There are a metric $d$ compatible with the topology
and positive constants $C$ and $\lambda$ with 
$\lambda >1$ such that for all $n>0$ and 
all points $x,y$ on a common edge of $X$, 
if $f^n$ maps the interval $[x,y]$ into an edge,
then $d(f^nx,f^ny)\geq C\lambda^n d(x,y)$.  
\item[Axiom 4.] ({\it Nonfolding})
$f^n|_{X-\mathcal{V}}$ is locally one-to-one
for every positive integer $n$.
\item[Axiom 5.]({\it Markov})
$f(\mathcal{V})\subseteq \mathcal{V}$.
\end{enumerate}

\begin{remarks}\label{2.1}
\begin{enumerate}
\item[(1)]
Axiom 0 means that
$X$ cannot be split into two nonempty,
closed, $f$-invariant subsets (\cite[\S1]{w2}).
\item[(2)]
We can define (without derivatives) an  arclength
from the assumed metric $d$ as follows.
Suppose $\gamma \colon [0,1]\to X$ where $\gamma$ 
is continuous and locally one-one. 
Let $t_1,\dots , t_{n-1}$ be all elements of 
$(0,1)$ which $\gamma$ maps to vertices of $X$. 
Define the length of $\gamma $ as 
$\sum_{i=1}^n d(x_{i-1},x_i)$. 
With this definition,
we can say that the  Expansion Axiom
means that there exists a metric
compatible with the topology of $X$ such that
there are constants $C>0$ and $\lambda >1$ such that
$f^n$ increases arclength
by a factor of at least $C\lambda^n$ 
(this was one formulation of Williams 
\cite{w1,w2}, except that his arclength 
tacitly was defined as usual with derivatives). 
Also, if we define a path metric $d'$ by 
setting   
$d'(x,y)$ to be the length of the shortest 
path from $x$ to $y$, then $d'$ is another 
metric compatible with the topology and 
still satisfies the Expansion Axiom. 
\end{enumerate}
\end{remarks}

\begin{standing assumption} \label{sa}
In this paper, we always assume that
$(X,f)$ satisfies Axioms 0 and 1.
\end{standing assumption}

For a given directed graph $X$ with a continuous map
$f\colon X\to X$,
let $\overline{X}$ be the inverse limit space
\[
\overline{X}
=X\overset{f}{\longleftarrow}
X\overset{f}{\longleftarrow}\cdots
=\bigl\{(x_0,x_1,x_2,\dots )\in
       \prod_0^\infty X \, |\, f(x_{n+1})=x_n \bigr\},
\]
and $\overline{f}\colon \overline{X}\to \overline{X}$
the induced homeomorphism defined by
\[
(x_0,x_1,x_2,\dots )\mapsto 
(f(x_0),f(x_1),f(x_2),\dots )=(f(x_0),x_0,x_1,\dots).
\]

Let $Y$ be a topological space and 
$g\colon Y\to Y$ a homeomorphism.
We call $Y$ a {\bf 1-dimensional generalized solenoid}
or {\bf $1$-solenoid} and $g$ a {\bf solenoid map}
if there exist a directed graph $X$
and a graph map $f\colon X\to X$
such that $(X,f)$ satisfies all six Axioms
and $(\overline{X},\overline{f})$
is topologically conjugate to $(Y,g)$.
We say that $(X,f)$ is a {\bf presentation}
of $(Y,g)$.
If $(X,f)$ satisfies all Axioms
except possibly the Flattening Axiom,
then we call $Y$ a {\bf branched solenoid}.

\begin{remarks} \label{2.2}
\begin{enumerate}
\item[(1)]
Williams  defined an $n$-dimensional
generalized solenoid $\overline{X}$
and a solenoid map $\overline{f}$
as the inverse limit of a system
$(X,f)$ satisfying Axioms 0-3 where $X$ is a
{\it differentiable $n$-dimensional branched manifold}
and $f\colon X\to X$ is
an {\it immersion} (\cite{w1, w2, w3}).
We generalize his systems in the topological category
for the $1$-dimensional case.
As a topological system,
every 1-solenoid of Williams is a
1-solenoid by our topological definition.
We will see the relation between
Williams' definition and the topological definition 
in Appendix B.
\item[(2)]
The Nonfolding Axiom  is the topological condition
we use in place of
the differentiable immersion condition.
\item[(3)]
If $(X,f)$ satisfies Axioms 0-4,
then there is a presentation $(X^\prime,f^\prime)$
satisfying Axioms 0-5 such that
$(\overline{X},\overline{f})$
is topologically conjugate to
$(\overline{X^\prime},\overline{f^\prime})$
(\cite[Proposition 3.5]{w1}). 
Williams proved \cite[3.5]{w1}
assuming the immersion condition,
but his proof goes through with our Axioms 0-4.
\end{enumerate}
\end{remarks}

\begin{example} \label{2.3}
Let $X$ be the unit circle on the complex plane.
Suppose that $1$ and $-1$ are the vertices of $X$,
and that the upper half circle $e_1$
and the lower half circle $e_2$
with counterclockwise direction are the edges of $X$.
Define $f,g:X\to X$ by
\[
f:z\mapsto z^2 \text{ and }
g:z\mapsto 
\begin{cases}
z^2 & \text{ if } z\in e_1\\ 
z^{-2} & \text{ if } z\in e_2.
\end{cases}
\]
Then $(X,f)$ satisfies all six Axioms, and
$(X,g)$ satisfies all Axioms
except the Nonfolding Axiom.
For $(X,g)$, $g^2$ is not locally one-to-one
at $\exp(\frac{\pi}{2}i)\in S^1$.
\end{example}

\begin{notation} \label{2.4}
Suppose that $(X,f)$ satisfies the Nonfolding Axiom
and the Markov Axiom,
$\mathcal{E}=\{e_1,\dots,e_n\}$ is the edge set
of $X$ with a given direction,
and $k$ is a positive integer.
For each edge $e_i\in \mathcal{E}$,
we can {\it give $e_i$ the partition
$\{I^{(k)}_{i,j}\}$, $1\le j \le j(i,k)$, for $f^k$}
such that
\begin{enumerate}
\item[(1)]
the initial point of $I^{(k)}_{i,1}$ is
the initial point of $e_i$,
\item[(2)]
the terminal point of $I^{(k)}_{i,j}$ is the initial
point of $I^{(k)}_{i,j+1}$ for $1\le j < j(i,k)$,
\item[(3)]
the terminal point of $I^{(k)}_{i,j(i,k)}$ is
the terminal point of $e_i$,
\item[(4)]
$f^k|_{\text{Int}{}I^{(k)}_{i,j}}$ is injective, and
\item[(5)]
$f^k(I^{(k)}_{i,j})={e_{i,j}^{(k)}}^{s(i,j,k)}$ where
$e_{i,j}^{(k)}\in \mathcal{E}$,
$s(i,j,k)=1$ if the direction of $f^k(I^{(k)}_{i,j})$
agree with that of $e_{i,j}^{(k)}$, and
$s(i,j,k)=-1$ if the direction of $f^k(I^{(k)}_{i,j})$
is reverse to that of $e_{i,j}^{(k)}$.
\end{enumerate}
\end{notation}

\begin{remarks}
\begin{itemize}
\item[(1)]
By the Nonfolding Axiom,
if $f^k(I^{(k)}_{i,j})=e_{i,j}^{\,\pm 1}$
for $1\le j < j(i,k)$,
then  $f^k(I^{(k)}_{i,j+1})=e_{i,j+1}$
cannot be $e_{i,j}^{\,\mp 1}$.
\item[(2)]
If $(X,f)$ satisfies all six axioms,
then there is a positive integer $m$ such that,
for every vertex $v$ of $X$ and 
every integer $k\ge m$,
there exist at most two edges
$e_{v,k,1}$ and $e_{v,k,2}$ such that,
for every
$I\in \{I_{i,j}^{(k)}\mid v\in f^k(I_{i,j}^{(k)}) \}$,
$f^{k}(I)=e_{v,k,l}^{\,\pm 1}$, $l=1$ or $2$.
\end{itemize}
\end{remarks}

\begin{definition} \label{2.5}
Suppose that
$(X,f)$ is a presentation of a branched solenoid,
that is, $(X,f)$ satisfies all Axioms
except possibly the Flattening Axiom,
and $\mathcal{E}$ is the edge set of $X$.
Then each edge $e_i\in \mathcal{E}$ has the partition
$\{ I^{(1)}_{i,j} \}$ for $f$, and
we can define an induced map
$\tilde{f}\colon \mathcal{E}\to \mathcal{E}^*$
by
\[
\tilde{f}\colon
e_i\mapsto e_{i,1}^{s(i,1)} e_{i,2}^{s(i,2)}
\cdots e_{i,j(i)}^{s(i,j(i))}
\]
where $e_{i,j}^{s(i,j)}=f(I^{(1)}_{i,j})$
for $1\le j \le j(i)$.
We call $\tilde{f}$ the {\bf substitution rule}
or the {\bf wrapping rule} associated to $f$.
\end{definition}

\begin{examples} \label{2.6}
Let $(X,f)$ and $(X,g)$ be given in Examples \ref{2.3}.
Then $\tilde{f}, \tilde{g}\colon
\mathcal{E}_X\to \mathcal{E}_X^*$ are given by
\[
\tilde{f}\colon e_1\mapsto e_1e_2,\,\,
          e_2\mapsto e_1e_2,
\text{ and }\,
\tilde{g}\colon e_1\mapsto e_1e_2,\,\,
     e_2\mapsto e_2^{-1}e_1^{-1}.
\]

To establish some notation,
we give Figure \ref{xcv1}
to represent
the presentations $(X,f)$ and $(X,g)$
with the wrapping rules $\tilde{f}$ and $\tilde{g}$,
respectively.
\begin{figure}[ht]
\centering
\begin{minipage}[c]{0.35\textwidth}
\centering \includegraphics[height=2.9cm]{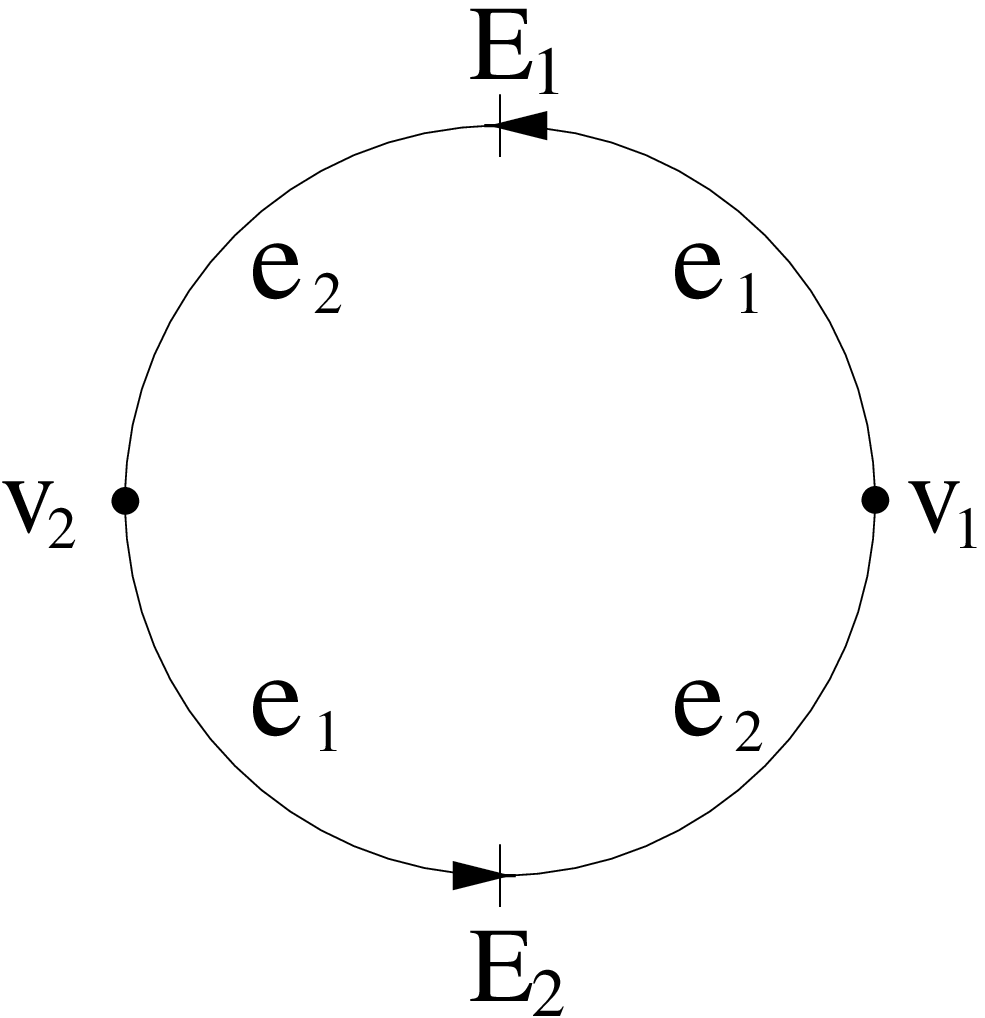}
\end{minipage}
\begin{minipage}[c]{0.35\textwidth}
\centering \includegraphics[height=2.9cm]{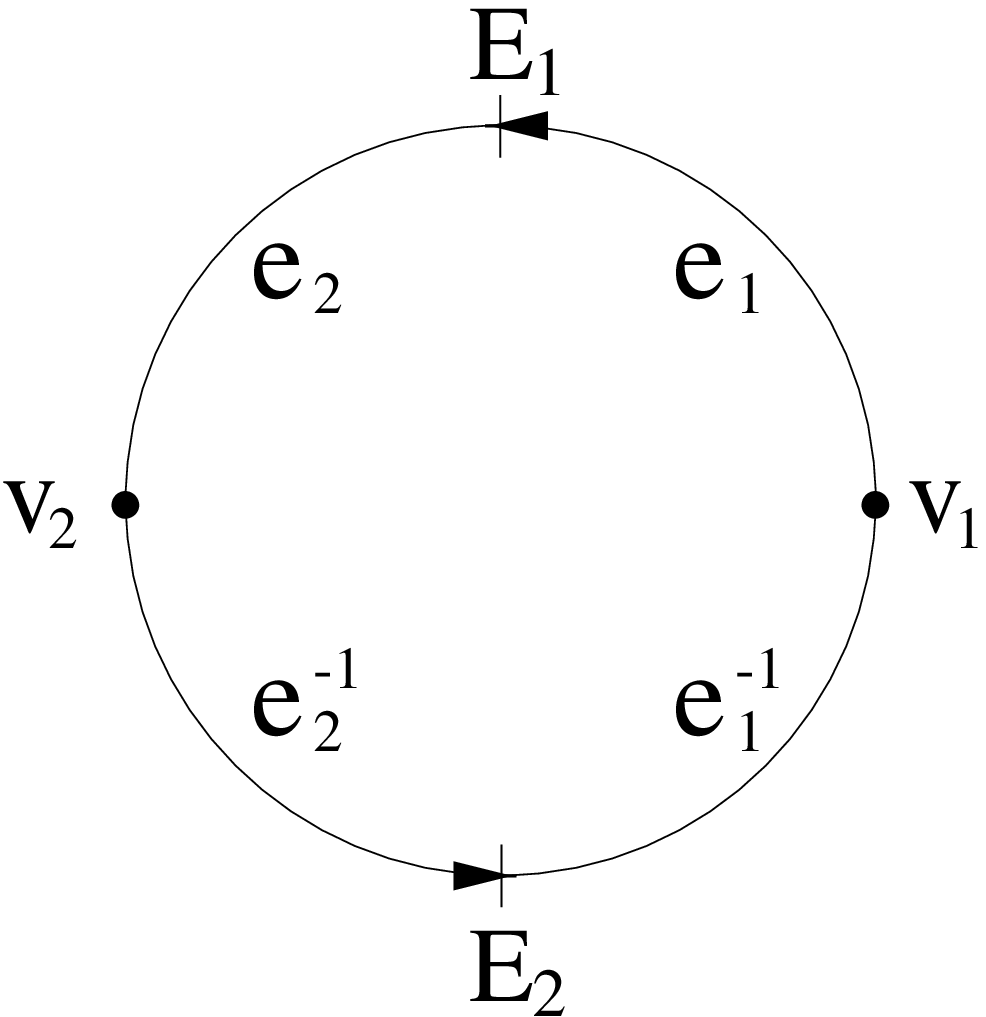}
\end{minipage}
\caption{$(X,f)$ and $(X,g)$
with the wrapping rules $\tilde{f}$ and $\tilde{g}$,
respectively.}
\label{xcv1}
\end{figure}

Similarly,
if $(Y,h)$ is given by Figure \ref{xcv2},
\begin{figure}[ht]
\centerline{\scalebox{.28}{\includegraphics{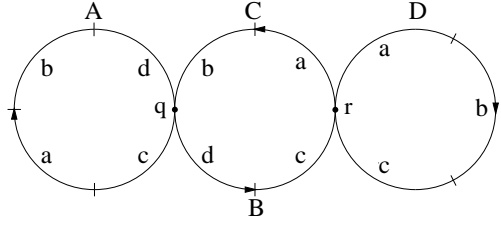}}}
\caption{$(Y,h)$ with wrapping rule $\tilde{h}$.}
\label{xcv2}
\end{figure}
then the wrapping rule 
$\tilde{h}\colon \mathcal{E}_Y\to \mathcal{E}_Y^*$
is given by
\[
a\mapsto cabd, \quad 
b\mapsto dc, \quad
c\mapsto ab, \quad
d\mapsto abc.
\]

Note that the two vertices $q,r$ of $Y$
have $h$-period $2$.
If $U_q$ and $U_r$ are
sufficiently small neighborhoods of $q$ and $r$,
respectively,
then $h^2(U_q)$ and $h(U_r)$ are intervals.
So $(Y,h)$ satisfies the Flattening Axiom.
\end{examples}

\begin{lemma} \label{2.7}
Suppose that $(X,f)$ satisfies Axioms $4$ and $5$.
Then there is a positive integer $l$ such that,
for each edge $e_i$ of $X$ and
every positive integer $m$,
if $\{I^{(l)}_{i,1},\dots , I^{(l)}_{i,j(i,l)}\}$ and
$\{I^{(lm)}_{i,1}, \dots ,I^{(lm)}_{i,j(i,lm)}\}$
are partitions of $e_i$ for $f^l$ and $f^{lm}$,
respectively,
then $f^l(I^{(l)}_{i,1})=f^{lm}(I^{(lm)}_{i,1})$
and $f^l(I^{(l)}_{i,j(i,l)})=
f^{lm}(I^{(lm)}_{i,j(i,lm)})$.
\end{lemma}
\begin{proof}
Since $\mathcal{V}$ is a finite set and
$f(\mathcal{V})\subset \mathcal{V}$,
every vertex of $X$ is eventually periodic, and
there is a positive integer $l_1$ such that
$f^{l_1}(v)=f^{l_1m}(v)$ 
for every $v\in \mathcal{V}$ and 
every positive integer $m$.

If $e_i$ is an edge of $X$ beginning at $v$,
and $\{ I^{(l_1)}_{i,1}, \dots ,
          I^{(l_1)}_{i,j(i,l_1)}\}$
is the partition of $e_i$ for $f^{l_1}$,
then $f^{l_1}(I^{(l_1)}_{i,1})=e_j^s$ such that
$f^{l_1}(v)\in e_j$.
Since $\mathcal{E}$ is a finite set
and $f^{l_1}(v)$ is a fixed point of $f^{l_1}$,
there is a positive integer $l_2$ such that,
for every positive integer $m$,
if $e_j$ has partitions
$\{I^{(l_2)}_{j,1},\dots\}$ for $f^{l_2}$ and 
$\{I^{(l_2m)}_{j,1},\dots\}$ for $f^{l_2m}$,
then
$f^{l_2}(I^{(l_2)}_{j,1})=f^{l_2m}(I^{(l_2m)}_{j,1})$.
This shows that 
$f^{l_1l_2}(I^{(l_1l_2)}_{i,1})
=f^{l_1l_2m}(I^{(l_1l_2m)}_{i,1})$
for every positive integer $m$.
By the same argument, we can choose a positive 
integer $l^\prime$ for $e_i$ such that
$f^{l^\prime}(I^{(l^\prime)}_{i,j(i,l^\prime)})=
f^{l^\prime m}(I^{(l^\prime m)}_{i,j(i,l^\prime m)})$.

Let $\ell_i$ be the least common multiple of
$l_1l_2$ and $l^\prime$ for each edge $e_i$, and
$l$  the least common multiple of these $\ell_i$'s.
Then we have
$f^l(I^{(l)}_{i,1})=f^{lm}(I^{(lm)}_{i,1})$ and
$f^l(I^{(l)}_{i,j(i,l)})=
     f^{lm}(I^{(lm)}_{i,j(i,lm)})$
for every positive integer $m$.
\end{proof}

\begin{lemma} \label{2.8}
Suppose that $(X,f)$ satisfies Axioms $3,4$ and $5$.
Then there exist  a positive integer $l$ and
$\epsilon >0$ such that, for all $x,y\in {X}$,
$f^l(x) \ne f^l(y)$ implies that
there is a nonnegative integer $K$ such that
$d({f}^K(x),{f}^K(y))\ge\epsilon$.
\end{lemma}
\begin{proof}
For convenience, we will take the metric $d$ on $X$
so that $d(x,y)$ is the length of the shortest path
between $x$ and $y$, as explained in Remarks \ref{2.1}.

Let $l$ be the integer given in Lemma~\ref{2.7}.
So each edge $e_i$ has the partition
$\{I^{(l)}_{i,j}\}$ for $f^l\colon X\to X$
as in Notation \ref{2.4}.
Without loss of generality, we suppose
each $j(i,l)\ge 3$.
Let $\mathcal{P}$ be the collection of
the intervals $I^{(l)}_{i,j}$.
First choose
$\epsilon_1>0$ so small that
\begin{itemize}
\item[(i)]
each interval $I^{(l)}_{i,j}$
has length larger than $2\epsilon_1>0$,
\item[(ii)]
if $f^l(x)=a\in I^{(l)}_{i,1}$,
$v$ is the initial point of $e_i$,
and $d(v,a)< \epsilon_1$,
then $f^l(a)\in I^{(l)}_{i,1}$, and
\item[(iii)]
if $f^l(y)=b\in I^{(l)}_{i,j(i,l)}$,
$w$ is the terminal point of $e_i$,
and $d(w,b)< \epsilon_1$,
then $f^l(b)\in I^{(l)}_{i,j(i,l)}$.
\end{itemize}
Then choose $\epsilon$ such that
$0<\epsilon<\epsilon_1$ and
for every $x$ in the compact set
$\bigcup\limits_{
\substack{ i\\ 1<j<j(i,l)}}I^{(l)}_{i,j}$
and every $y\in X$,
if $0<d(x,y)<\epsilon$, then
\begin{itemize}
\item[(iv)]
$0< d(f^l(x),f^l(y))<\epsilon_1$ and
\item[(v)]
the interval $[f^l(x),f^l(y)]$ contains at most
one vertex.
\end{itemize}
Note that $f^l(x)\ne f^l(y)$ comes from
the Nonfolding Axiom.

If $f^l(x)=a\ne b=f^l(y)$ and $d(a,b)< \epsilon_1$,
then $a$ and $b$ lie on the same or adjacent
elements of $\mathcal{P}$.
So we have two cases:
\begin{enumerate}
\item[(1)]
$a$ and $b$ are end points of an interval
of length less than $\epsilon_1$
containing a vertex $v$, or
\item[(2)]
the interval $[a,b]$ of length less than $\epsilon_1$
does not contain any vertex of $X$.
\end{enumerate}

For case (1), by the condition (i),
$d(a,b)< \epsilon_1$ implies that
$a$ and $b$ are contained in two different intervals
among $I^{(l)}_{i,1}$, $I^{(l)}_{i,j(i,l)}$,
$I^{(l)}_{n,1}$, and $I^{(l)}_{n,j(n,l)}$.
For brevity, let's assume
$a\in I_{i,1}$ and $b\in I_{n,1}$.
Then by the condition (ii),
$f^l(a)\in I^{(l)}_{i,1}$ and
$f^l(b)\in I^{(l)}_{n,1}$.
If $v$ is the vertex of $X$ contained in 
$I^{(l)}_{i,1}\cap I^{(l)}_{n,1}$,
then $f^l$ maps $[v,a]$ into $I_{i,1}$.
So
\[
d(f^l(a),v) \ge c\lambda^l \cdot d(a,v)
\]
where $c$ and $\lambda$ are
the expansion constants.
Similarly 
$d(f^l(b),v) \ge c\lambda^l \cdot d(b,v)$. 
Let $k$ be the smallest positive integer such that
$d(v,f^{lk}(a))\ge \epsilon_1$ or
$d(v,f^{lk}(b))\ge \epsilon_1$.
Then by induction using (ii) and (iii),
we have for $0<s\le k$ that
$f^{sl}$ sends $[v,a]$ injectively into $I_{i,1}$
and $[v,b]$ injectively into $I_{n,1}$.
Therefore we have
$d(f^{kl}(a),f^{kl}(b))\ge\epsilon_1>\epsilon$.

For case (2),
let $k$ be the smallest positive integer such that
$[f^{kl}(a),f^{kl}(b)]$ contains a vertex.
It follows from the Nonfolding Axiom
that $f^{(k-1)l}(a)\ne f^{(k-1)l}(b)$.
If there exists $i$, $0<i<k$, such that
$d(f^{il}(a), f^{il}(b))>\epsilon$,
then we are done, so suppose not.
Then $f^{kl}(a)$ and $f^{kl}(b)$ are
endpoints of an interval of length less than
$\epsilon_1$ containing a vertex.
Hence we may conclude the proof
by applying the argument of case (1).
\end{proof}

\begin{definition}[{\cite[\S3.5]{ro}}] \label{2.9}
A homeomorphism $h$ on a metric space $Y$
is {\bf expansive} if there is an $\epsilon> 0$ 
such that, for all $x\ne y\in Y$,
there is an integer $n$ such that
\[
d(h^n(x),h^n(y))\ge \epsilon.
\]
\end{definition}	

For a solenoid $\overline{X}$ presented by $(X,f)$,
define a metric $\bar{d}$ on $\overline{X}$ by
\[
\bar{d}(x,y)=\sum\limits_{i=0}^{\infty}
 \frac{d(x_i,y_i)}{2^i}
\]
where $x=(x_0,x_1,\dots), y=(y_0,y_1,\dots)
\in \overline{X}$, and
$d$ is a metric on $X$ compatible with
the topology of $X$ such that
$f\colon X\to X$ is an expansion with respect to
$d$ as in Remarks \ref{2.1}.(2).

\begin{proposition} \label{2.10}
If $(X,f)$ satisfies Axioms $3,4$ and $5$,
then $\overline{f}\colon \overline{X}\to \overline{X}$
is expansive.
\end{proposition}
\begin{proof}
For a pair of points
$x=(x_0,x_1,\dots)\ne y=(y_0,y_1,\dots)
\in \overline{X}$,
there is a nonnegative integer $N$ such that
$x_n\ne y_n$ for all $n\ge N$.

Let $l$ and $\epsilon > 0$ be given in Lemma~\ref{2.8}.
Then $x_{N+l}\ne y_{N+l}$ and
$f^l(x_{N+l})\ne f^l(y_{N+l})$
imply that there exists $K\ge 0$ such that
$d\bigl({f}^K(x_{N+l}),{f}^K(y_{N+l})\bigr)\ge\epsilon$
by Lemma~\ref{2.8}.
Therefore we have 
\[
\bar{d}\bigl(\overline{f}^{K-N-l}(x),
   \overline{f}^{K-N-l}(y) \bigr) 
= d\bigl({f}^K(x_{N+l}),{f}^K(y_{N+l})\bigr)
+\sum\limits_{i > K-N-l}
\frac{d\bigl(x_{i},y_{i})\bigr)}{2^{N+l-K+i}}
>\epsilon,
\]
and this proves that
$\overline{f}\colon \overline{X}\to \overline{X}$
is expansive.
\end{proof}

\begin{example} \label{2.11}
We need the Nonfolding Axiom in Proposition~\ref{2.10}.

For $f,g \colon S^1\to S^1$ in Examples \ref{2.3},
it follows from Proposition~\ref{2.10}
that $\overline{f}$ is expansive.
For $\overline{g}$ with given $\epsilon > 0$
and $x=(x_0,x_1,\dots), y=(y_0,y_1,\dots)
\in (\overline{S^1},\overline{g})$, let
\[
x_n=\exp (\frac{\pi}{2^n}i +
          \frac{1}{2^{n+2}}\epsilon i)
\text{ and }
y_n=\exp (\frac{\pi}{2^n}i -
          \frac{1}{2^{n+2}}\epsilon i).
\]
Then, for a natural Riemannian metric $d$ on $S^1$,
\[\bar{d}(\overline{g}^k(x),\overline{g}^k(y))
=\frac{1}{2^{|k|}} d(x,y)=\frac{1}{2^{|k|}}\epsilon
\]
for every integer $k$,
and $\overline{g}$ is not expansive.
\end{example}

\subsection*{SFT covers}
We will review the standard construction of
a shift of finite type (SFT) cover for
$1$-dimensional branched solenoids.

Suppose that
$(X,f)$ is a presentation of a branched solenoid,
and $\mathcal{E}=\{e_1,\dots,e_n\}$
is the edge set of $X$.
Let $\{ I_{i,j}\mid 1\le i \le n, 1\le j \le j(i) \}$
be the partition of $\mathcal{E}$ for $f$,
and 
$\tilde{f}\colon \mathcal{E} \to \mathcal{E}^*$
the wrapping rule associated to $f$ given by
\begin{equation} \label{e1}
\tilde{f}\colon
e_i=I_{i,1}\cdots I_{i,j(i)}
\mapsto e_{i,1}^{s(i,1)} e_{i,2}^{s(i,2)}
         \cdots e_{i,j(i)}^{s(i,j(i))}
\end{equation}
where $e_{i,j}^{s(i,j)}=f(I_{i,j})$
and $s(i,j)=\pm 1$ denotes the direction.
The {\it adjacency matrix} $M$ of
$(\mathcal{E},\tilde{f})$ is given by
\[
M(i,j)=\# \{ I_{i,l}\mid f(I_{i,l})=e_j^{\pm 1}\}.
\]

We may view $M$ as the adjacency matrix of
a directed graph whose vertex set is $\mathcal{E}$
and whose edge set is
$\mathcal{A}=\{ I_{i,j}\mid 1\le i\le n, 
                           1\le j \le j(i) \}$,
the partition of $\mathcal{E}$ for $f$,
where outgoing edges from the `vertex' $e_i$
are named $I_{i,1},\dots, I_{i,j(i)}$.

Now we can give $(\overline{X},\overline{f})$
a two-sided SFT cover $(\Sigma_X,\sigma_X)$
defined from the {\sl alphabet} $\mathcal{A}$
and the adjacency matrix $M_X$.
The shift space $\Sigma_X$ is the subset of
$\mathcal{A}^{\mathbb{Z}}$ whose forbidden blocks are
$\{ I_{i,j}I_{k,l}\mid
I_{k,l}\nsubseteq {f}(I_{i,j})=e_{i,j}^{s(i,j)} \}$
from the formula $(1)$.
Therefore $\Sigma_X$ is a 1-step subshift of
finite type,
and a word $w=I_{a(0)}I_{a(1)}\cdots I_{a(l)}$ is
{\sl allowed} in $\Sigma_X$ if and only if
$\bigcap\limits_{j=0}^l {f}^{-j}(I_{a(j)})$ has
nonempty interior.

For each point
$w=\cdots I_{a(-1)}I_{a(0)}I_{a(1)}\cdots
\in \Sigma_X$
and the canonical projection map onto
the zeroth coordinate
$\pi\colon \overline{X}\to X$, 
there is a unique corresponding point
\[
x_w=(x_0,x_1,\dots)=\bigcap_{j=-\infty}^\infty 
\overline{{\bar{f}}^{-j}
 \bigl(\pi^{-1}(I_{a(j)})\bigr)} \in \overline{X}
\]
such that 
$x_i\in I_{a(-i)}$ and $f^i(x_0)\in I_{a(i)}$
for $i\ge 0$.
Hence there is a well-defined semiconjugacy
$p\colon\Sigma_X\to\overline{X}$
defined by $w\mapsto x_w$.
It is not difficult to check
${\overline{f}}\circ p=p \circ \sigma_X$
(\cite[\S6.5]{lm}).

\begin{proposition}[{\cite[\S3.D]{b2}}]\label{2.12}
Suppose that $(X,f)$ satisfies all six Axioms
except possibly the Flattening Axiom. 
Let $p$ and $\Sigma_X$ be as above. Then
\begin{enumerate}
\item[(1)]
$p\colon\Sigma_X\to\overline{X}$ is finite-to-one.
\item[(2)]
$p$ is one-to-one on
$\Sigma_X \backslash\bigcup\limits_{m=0}^\infty 
p^{-1}\circ \pi^{-1}\circ f^{-m-1}({\mathcal{V}})$
where $\pi\colon \overline{X}\to X$
is the projection map to the zeroth coordinate space
and $\mathcal{V}$ is the vertex set of $X$.
\item[(3)]
$(\overline{X},\overline{f})$ and
$(\Sigma_X,\sigma_X)$ have the same entropy.
\end{enumerate}
\end{proposition}

\begin{lemma}[{\cite[1.6]{w2}}]\label{2.19}
If $(X,f)$ satisfies all six Axioms,
and $I$ is an interval in $X$,
then $X\subset f^m(I)$ for some $m\ge 0$.
\end{lemma}

We remark that the proof of \cite[1.6]{w2}
still works in our topological setting.
Then for all $I_{i,j}, I_{k,l}\in \mathcal{A}$,
$I_{k,l}\subset f^m(I_{i,j})$
for some positive integer $m$,
and we have the following proposition.

\begin{proposition}\label{2.17}
If $(X,f)$ satisfies all six Axioms,
then $(\Sigma_X,\sigma_X)$ is a mixing SFT.
\end{proposition}

\begin{examples} \label{2.13}
Let $(X,f)$ and $(Y,h)$ be as in Example~\ref{2.6}.
Recall that
$\tilde{f}\colon {\mathcal{E}}_X\to {\mathcal{E}}_X^*$
and
$\tilde{h}\colon {\mathcal{E}}_Y\to {\mathcal{E}}_Y^*$
are given by
\begin{align*}
\tilde{f}&\colon
e_1\mapsto e_1e_2, \text{ }
e_2\mapsto e_1e_2,\\
\tilde{h}&\colon
a\mapsto cabd,\text{ } b\mapsto dc,\text{ and }
c\mapsto abc,\text{ } d\mapsto ab.
\end{align*}

The SFT covers of $(\overline{X},\overline{f})$
and $(\overline{Y},\overline{h})$
are given by the following matrices
\[
M_X=
\begin{pmatrix}
1&1\\ 
1&1
\end{pmatrix}
\text{ and }
M_Y=
\begin{pmatrix}
1&1&1&1\\
0&0&1&1\\
1&1&0&0\\
1&1&1&0
\end{pmatrix}.
\]
\end{examples}

\begin{example}\label{2.18}
We need the Flattening Axiom
in Lemma \ref{2.19} and Proposition \ref{2.17}.
Let $Z$ be a wedge product of two circles $a$ and $b$,
and $g\colon Z\to Z$ given by
$a\mapsto aa$ and $b\mapsto bb$.
Then $(Z,g)$ does not satisfy the Flattening Axiom,
and the adjacency matrix $M_Z$
for $(\mathcal{E},\tilde{g})$
is 
$M_Z=\begin{pmatrix}
2&0\\0&2
\end{pmatrix}$.
So $(\Sigma_Z,\sigma_Z)$ is not irreducible.
\end{example}

\section{Shift equivalence}\label{S3}

We define shift equivalence of directed graphs
with graph maps,
and show that
conjugacy of branched solenoids is equivalent 
to shift equivalence of their presentations
and that certain conjugacies of solenoids
lift uniquely to conjugacies of associated SFT covers.

\begin{definition} \label{3.1}
Suppose that $X$ and $Y$ are directed graphs, and
that $f\colon X\to X$ and
$g\colon Y\to Y$ are graph maps.
We say that $f$ and $g$ are
{\bf shift equivalent with lag${}\, m$}
if there are continuous maps
$r\colon X\to Y$, $s\colon Y\to X$
and a positive integer $m$ 
such that 
\[
r\circ f=g\circ r,\quad f\circ s=s\circ g,\quad
f^m=s\circ r,\quad g^m=r\circ s.
\]
\end{definition}

If $(X,f)$ and $(Y,g)$ are
shift equivalent with lag${}\,m$ by continuous maps
$r\colon X\to Y$ and $s\colon Y\to X$,
then define induced maps
$\overline{r}\colon \overline{X}\to \overline{Y}$ 
and 
$\overline{s}\colon \overline{Y}\to \overline{X}$ by
\begin{align*}
\overline{r}&\colon
(x_0,\dots,x_m,x_{m+1},\dots)\mapsto
(r(x_m),r(x_{m+1}),\dots)\\
\overline{s}&\colon
(y_0,y_1,y_2,\dots)\mapsto
(s(y_0),s(y_1),s(y_2),\dots).
\end{align*}
We can easily check that
$\overline{r}$ and $\overline{s}$
are topological conjugacies of
$\overline{f}$ and $\overline{g}$ such that
$\overline{s}\circ \overline{r}=Id$ on ${\overline{X}}$
and $\overline{r}\circ \overline{s}=Id$ on
${\overline{Y}}$.

\begin{remark} \label{3.2}
It is possible that the shift equivalence map
$r\colon X\to Y$ is not a graph map,
that is, a vertex of $X$ may not be
mapped to a vertex of $Y$.
\end{remark}

\begin{lemma} \label{3.3}
Suppose that
$(X,f)$ and $(Y,g)$ satisfy the Markov Axiom,
and they are shift equivalent to each other
with lag$\,{}m$ by continuous maps
$r\colon X\to Y$ and $s\colon Y\to X$. 
Then the vertex sets
${\mathcal{V}}_X$ and ${\mathcal{V}}_Y$
of $X$ and $Y$, respectively,
can be enlarged to $\mathcal{V}_X^\prime$
and $\mathcal{V}_Y^\prime$, respectively,
so that 
$$
f(\mathcal{V}_X^\prime)\subset
\mathcal{V}_X^\prime,\text{ }
g(\mathcal{V}_Y^\prime)\subset
\mathcal{V}_Y^\prime,\text{ }
r(\mathcal{V}_X^\prime)\subset
\mathcal{V}_Y^\prime,\text{ and }
s(\mathcal{V}_Y^\prime)\subset
\mathcal{V}_X^\prime.
$$
\end{lemma}
\begin{proof}
For the vertex sets
${\mathcal{V}}_X$ and ${\mathcal{V}}_Y$,
let $\mathcal{V}_Y^\prime=
\mathcal{V}_Y\cup r(\mathcal{V}_X)$
be the set of enlarged vertices in $Y$ and
$\mathcal{V}_X^\prime=
\mathcal{V}_X\cup s(\mathcal{V}_Y^\prime)$.
Then
\begin{align*}
g(\mathcal{V}_Y^\prime)&\subset 
g(\mathcal{V}_Y)\cup g\circ r(\mathcal{V}_X)\subset
\mathcal{V}_Y\cup r\circ f(\mathcal{V}_X)\subset
\mathcal{V}_Y\cup r(\mathcal{V}_X)=
\mathcal{V}_Y^\prime,\\
f(\mathcal{V}_X^\prime)&\subset
f(\mathcal{V}_X)\cup f\circ s(\mathcal{V}_Y^\prime)
\subset
\mathcal{V}_X\cup s\circ g(\mathcal{V}_Y^\prime)\subset
\mathcal{V}_X\cup s(\mathcal{V}_Y^\prime)
=\mathcal{V}_X^\prime,\\
r(\mathcal{V}_X^\prime)
&=r(\mathcal{V}_X)\cup r\circ s(\mathcal{V}_Y^\prime)
=r(\mathcal{V}_X)\cup g^m(\mathcal{V}_Y^\prime)
\subseteq\mathcal{V}_Y^\prime,\text{ and}\\
s(\mathcal{V}_Y^\prime)
&=s(\mathcal{V}_Y)\cup f^m(\mathcal{V}_X)
\subseteq\mathcal{V}_X^\prime
\end{align*}
prove the Lemma.
\end{proof}

\begin{lemma}[Ladder Lemma]\label{2.14}
Suppose that $(X,f)$ and $(Y,g)$ satisfy all
Axioms except possibly the Flattening Axiom,
and that $\phi\colon \overline{X}\to \overline{Y}$
is a continuous map such that
$\overline{g}\circ \phi=\phi\circ\overline{f}$.
Then there is a continuous map $r\colon X\to Y$
and a nonnegative  integer $n$ such that
$g\circ r=r\circ f$ and 
$\phi(x_0,x_1,\dots)=(r(x_n),r(x_{n+1}),\dots)$.
\end{lemma}

\begin{remark} \label{2.15}
Williams (\cite[\S3]{w2}) proved the Ladder Lemma
under the hypotheses
that $f$ and $g$ are nonwandering expansions
which are immersions of
differentiable branched $1$-manifolds
and  satisfy the Flattening Axiom.
Our assumptions are weaker than his conditions
as we don't need the Flattening Axiom,
but the ideas of the proof are
essentially those given by Williams.
The additional complication of our setting 
is addressed in Lemma \ref{2.8}.
\end{remark}

Let $X_i$ and $Y_i$ be the $i$th coordinate
spaces of $\overline{X}$ and $\overline{Y}$,
respectively,
and $\pi_i$ the projection map
from the branched solenoids onto
their $i$th coordinate spaces.

\begin{lemma} \label{2.16}
There is a positive integer
$n$ such that, for $a,b\in \overline{X}$,
if $\pi_n(a) = \pi_n(b)$, 
then $\pi_0\circ \phi(a) = \pi_0\circ \phi(b)$.
\end{lemma}
\begin{proof}
By Lemma \ref{2.8}, we can choose $\epsilon> 0$
and $l\in\mathbb{N}$ such that, for all $x,y\in Y$,
if $g^l(x)\ne g^l(y)$,
then there exists a nonnegative integer $K$
such that \[d(g^K(x),g^K(y))\ge \epsilon.\]
Choose a  nonnegative integer $k$ and
$\delta > 0$ such that,
for $a=(a_0,a_1,\dots)$ and
$b=(b_0,b_1,\dots)\in \overline{X}$, 
$a_k=b_k$ implies $d(a,b)< \delta$, and
$d(a,b)< \delta$ implies $d(\phi(a),\phi(b))<\epsilon$.

Now suppose $a_k=b_k$.
From $d(\overline{f}^m(a),\overline{f}^m(b))\le
d(a,b)<\delta$ for every nonnegative integer $m$,
we have
\[
d(\phi\circ\overline{f}^m(a),
            \phi\circ\overline{f}^m(b))
=d(\overline{g}^m\circ\phi(a),
                \overline{g}^m\circ\phi(b))
<\epsilon,
\]
and if $x_0=\pi_0\circ \phi(a)$ and
$y_0=\pi_0\circ \phi(b)$, then
for every nonnegative integer $m$
\[
d({g}^m(x_0),{g}^m(y_0))<\epsilon.
\]
Therefore we have that $a_k=b_k$ implies
$g^l\circ\pi_0\circ\phi(a)=g^l\circ\pi_0\circ\phi(b)$,
and equivalently $a_{k+l}=b_{k+l}$ implies
$\pi_0\circ\phi(a)=\pi_0\circ\phi(b)$.
\end{proof}

\begin{proof}[Proof of Ladder Lemma]
Let $n$ be as in Lemma \ref{2.16}.
Define $r_k\colon X\to Y$ by
$r_k=\pi_k \circ \phi \circ \pi_{n+k}^{-1}$.
Then $r_k$ is well-defined by our choice of $n$.
Now show that $r_k=r_0$ for every positive integer $k$.
For $x=(x_0,\dots,x_n,\dots)\in\overline{X}$,
\begin{align*}
r_k(x_n)
&=\pi_k\circ\phi\circ\pi_{n+k}^{-1}(x_n)
=\pi_k\circ\phi\circ\overline{f}^k(x)\\
&=\pi_k\circ\overline{g}^k\circ\phi(x)
=\pi_0\circ\phi(x)
=\pi_0\circ\phi\circ\pi_n^{-1}(x_n)\\
&=r_0(x_n).
\end{align*}

To show the continuity of $r=r_0$,
let $\delta>0$ and $\epsilon>0$ be
as in Lemma \ref{2.16},
and choose $\delta^\prime>0$ such that
if $a_n, b_n$ are elements in $X$
with $d(a_n,b_n)<\delta^\prime$,
then there exist $a,b\in \overline{X}$ 
with $\pi_n(a)=a_n$ and $\pi_n(b)=b_n$ such that
$d(a,b)<\delta$.
Then we have
\[
d(a_n,b_n)<\delta^\prime
\implies
d(\pi_0\circ\phi(a),\pi_0\circ\phi(b))
=d(r(a_n),r(b_n)) < \epsilon,
\]
and $r\colon X\to Y$ is continuous.
That $\phi(x_0,\dots,x_n,x_{n+1},\dots)
=(r(x_n),r(x_{n+1}),\dots)$ is trivial
by the construction of $r\colon X\to Y$.
\end{proof}

The proof of the following proposition is
easy from the Ladder Lemma.
So we omit the proof.

\begin{proposition}[{\cite[Theorem 3.3]{w2}}]\label{3.4}
Suppose that $(X,f)$ and $(Y,g)$ satisfy
all Axioms except possibly the Flattening Axiom.
Then $\phi\colon\overline{X}\to \overline{Y}$ is
a topological conjugacy if and only if
there exists a shift equivalence
$(r,s)$ of $f$ and $g$ such that
$\phi=\overline{r}$.
\end{proposition}

Let $(\Sigma_X,\sigma_X)$ and $(\Sigma_Y,\sigma_Y)$
be the SFT covers of 
$(\overline{X},\overline{f})$ and 
$(\overline{Y},\overline{g})$
defined by $(\mathcal{E}_X,\tilde{f})$
and $(\mathcal{E}_Y,\tilde{g})$, respectively,
as in \S\ref{S2}.

\begin{theorem} \label{3.5}
Suppose that
$(X,f)$ and $(Y,g)$ satisfy Axioms $3,4,5$, and
they are shift equivalent to each other
with $\text{lag}\,{}m$ by graph maps $r$ and $s$.
Then the conjugacy
$\overline{r}\colon \overline{X}\to \overline{Y}$
lifts to a unique conjugacy $\tilde{r}$ of
$(\Sigma_X,\sigma_X)$ and $(\Sigma_Y,\sigma_Y)$.
\end{theorem}
\begin{proof}
We will define a sliding block code
$\phi_{r}\colon \Sigma_X \to \Sigma_Y$ induced by $r$, 
and show that $\phi_{r}$ is the required conjugacy
$\tilde{r}$.

Let $\mathcal{E}_X$ and $\mathcal{E}_Y$
denote the edge sets of $X$ and $Y$, respectively,
$\mathcal{A}_X=\{I_{i,j}\}$ and
$\mathcal{A}_Y=\{J_{k,l}\}$
the alphabets of $(\Sigma_X,\sigma_X)$ and
$(\Sigma_Y,\sigma_Y)$, respectively,
and $p_X\colon \Sigma_X\to \overline{X}$
and $p_Y\colon \Sigma_Y\to \overline{Y}$
the semiconjugacies.
Then $I_{i,j}$ is a subset of $e_i\in \mathcal{E}_X$
such that $f(I_{i,j})=e_{i,j}\in \mathcal{E}_X$
(ignoring the direction).
Note that if
$a=\cdots I_{a(-1)}I_{a(0)}I_{a(1)}\cdots
\in \Sigma_X$, $x_a=(x_0,x_1,\dots)=p_X(a)
\in \overline{X}$, and $I_{a(i)}\subset e_{a_{i}}$,
then $x_i\in I_{a(-i)}$ and $f^i(x_0)\in I_{a(i)}$
for every nonnegative integer $i$.

Let $\mathcal{C}_X=\{ a\in \Sigma_X \mid 
p_X^{-1}\circ p_X(a)=\{ a \} \}$
and $\mathcal{C}_Y=\{ b\in \Sigma_Y \mid 
p_Y^{-1}\circ p_Y(b)=\{ b \} \}$.
Then by Proposition \ref{2.12}, 
$\mathcal{C}_X$ and $\mathcal{C}_Y$ are
dense in $\Sigma_X$ and $\Sigma_Y$, respectively.

\subsubsection*{Step 1}
Show that
$\overline{r}\circ p_X(\mathcal{C}_X)=
p_Y(\mathcal{C}_Y)$.

By Proposition \ref{2.12}, 
$a\in \mathcal{C}_X$ if and only if
$\pi_0\circ p_X(a)\notin
\bigcup\limits_{n=0}^{\infty}f^{-n-1}(\mathcal{V}_X)$
where $\pi_i$ is the projection map
from the branched solenoids
to their $i$th coordinate spaces. 
So we have $\overline{f}^m\circ p_X(\mathcal{C}_X)
=p_X(\mathcal{C}_X)$.

If $\overline{r}\circ p_X(a)\notin
p_Y(\mathcal{C}_Y)$ for some $a\in \mathcal{C}_X$,
then $r\circ\pi_m\circ p_X(a) \in
\bigcup\limits_{n=0}^{\infty}g^{-n-1}(\mathcal{V}_Y)$
and
\[
g^{n+1}\circ r\circ\pi_m\circ p_X(a)
=r\circ f^{n+1}\circ\pi_m\circ p_X(a)\in \mathcal{V}_Y
\]
for some $n\ge 0$.
Since the shift equivalence maps
$r$ and $s$ are graph maps by Lemma~\ref{3.3}
and $s\circ r=f^m$, we have
\[
s\circ r\circ f^{n+1}\circ\pi_m\circ p_X(a)
=f^{n+m+1}\circ\pi_m\circ p_X(a)
=f^{n+1}\circ\pi_0\circ p_X(a)
\in \mathcal{V}_X,
\]
a contradiction.
Hence we have
$\overline{r}\circ p_X(\mathcal{C}_X)\subset
p_Y(\mathcal{C}_Y)$.
By the same argument, we can show that
$\overline{s}\circ p_Y(\mathcal{C}_Y)\subset
p_X(\mathcal{C}_X)$.
Then $\overline{r}\circ\overline{s}=\overline{g}^m$
and $\overline{g}^m(\mathcal{C}_Y)=\mathcal{C}_Y$
imply that
\[
p_Y(\mathcal{C}_Y)=
\overline{r}\circ\overline{s}\circ p_Y(\mathcal{C}_Y)
\subset
\overline{r}\circ p_X(\mathcal{C}_X)
\]
Therefore we have
$\overline{r}\circ p_X(\mathcal{C}_X)
=p_Y(\mathcal{C}_Y)$.

Now we have a well-defined bijective map
$(p_Y|_{\mathcal{C}_Y})^{-1}\circ\overline{r}\circ
p_X|_{\mathcal{C}_X} \colon
\mathcal{C}_X\to \mathcal{C}_Y$.
This map will define $\tilde{r}$ on ${\mathcal{C}_X}$.

\subsubsection*{Step 2}
Find a block map
$\Phi_r\colon
\mathcal{B}_{m+1}(\mathcal{C}_X)\to \mathcal{A}_Y$
where $\mathcal{B}_{m+1}(\mathcal{C}_X)$ 
is the set of all $(m+1)$-blocks in $\mathcal{C}_X$ 
such that for every
$a=\cdots I_{a(-m+i)}\cdots I_{a(-1+i)}I_{a(i)}\cdots
\in \mathcal{C}_X$
$(p_Y|_{\mathcal{C}_Y})^{-1}\circ\overline{r}\circ
p_X|_{\mathcal{C}_X}(a)_i
=\Phi_r(I_{a(-m+i)}\cdots I_{a(i)})$.

For $a=\cdots I_{a(-m)}\cdots I_{a(-1)}I_{a(0)}\cdots
\in \mathcal{C}_X$,
let $x_a=(x_0,x_1,\dots)=p_X(a)\in \overline{X}$ and
$y_a=(y_0,y_1,\dots)=\overline{r}(x_a)
\in p_Y(\mathcal{C}_Y)$.
Then $x_i\in I_{a(-i)}\subset e_{a(-i)}$,
$y_i=r(x_{i+m})$,
and there exists a unique
$\alpha=\cdots J_{\alpha(-1)} J_{\alpha(0)}\cdots
\in \mathcal{C}_Y$ such that
\[
(p_Y|_{\mathcal{C}_Y})^{-1}\circ\overline{r}\circ
p_X|_{\mathcal{C}_X}(a)=\alpha \text{ and }
p_Y(\alpha)=y_a.
\]

Let $\{I_{e,j}^{(k)}\}$ be the partition of
$e=e_{a(-m)}\in \mathcal{E}_X$ for
$f^k$, $1\le k \le m+1$.
Then each $I_{e,j}^{(k)}$ is contained in
a unique $I_{e,j^\prime}^{(k-1)}$,
and we have a unique descending sequence
\[
e\supset I_{a(-m)}=I_{e,u(1)}^{(1)}\supset
I_{e,u(2)}^{(2)}\supset\cdots \supset
I_{e,u(m+1)}^{(m+1)} 
\]
such that
$f^k(I_{e,u(k)}^{(k)})=e_{a(-m+k)}\in \mathcal{E}_X$,
$f^k(I_{e,u(k+1)}^{(k+1)})
=I_{a(-m+k)}\subset e_{a(-m+k)}$,
and $x_m\in I_{e,u(m+1)}^{(m+1)}$.

Since the shift equivalence maps
$r$ and $s$ are graph maps and $s\circ r=f^m$,
$f^m(I_{e,u(m)}^{(m)})=e_{a(0)}$
and $f^{m+1}(I_{e,u(m+1)}^{(m+1)})=e_{a(1)}$
imply that
$r(I_{e,u(m)}^{(m)})$ is contained in
a unique edge $\epsilon$ in $Y$
and that $r(I_{e,u(m+1)}^{(m+1)})$ is contained in
a unique path $J$  such that
$J\subset \epsilon$ and $g(J)\in \mathcal{E}_Y$.
Define a block map
$\Phi_r\colon
\mathcal{B}_{m+1}(\mathcal{C}_X)\to \mathcal{A}_Y$
by
\[
I_{a(-m)}\cdots I_{a(-1)} I_{a(0)}\mapsto J.
\]
Then the sliding block code ${\phi_r}$
induced by $\Phi_r$ maps $a$ to
$\beta=\cdots J_{\beta(0)} J_{\beta(1)}\cdots$
with $J_{\beta(i)}
=\Phi_r(I_{a(-m+i)}\cdots I_{a(i)})$.

To prove that ${\phi_r}=
(p_Y|_{\mathcal{C}_Y})^{-1}\circ\overline{r}\circ
p_X|_{\mathcal{C}_X}$,
we need only show that $p_Y(\beta)=y_a$,
that is, $y_i=r(x_{i+m})\in J_{\beta(-i)}$.
From the descending sequence for $x_{i+m}$
\begin{equation} \label{e2}
e_{a(-i-m)}\supset I_{a(-i-m)}\supset
I_{a(-i-m),u(2)}^{(2)}\supset\cdots \supset
I_{a(-i-m),u(m+1)}^{(m+1)}\ni x_{i+m},
\end{equation}
$J_{\beta(-i)}$ is the unique path in $Y$ such that
$r(I_{a(-i-m),u(m+1)}^{(m+1)})\subset J_{\beta(-i)}$
and $g(J_{\beta(-i)}) \in \mathcal{E}_Y$.
So $y_i=r(x_{i+m})$ is contained in 
$J_{\beta(-i)}$, and $p_Y(\beta)=y_a$.
Therefore
${\phi_r}(a)=\beta=\alpha=
(p_Y|_{\mathcal{C}_Y})^{-1}\circ\overline{r}\circ
p_X|_{\mathcal{C}_X}(a)$
by the definition of $\mathcal{C}_Y$,
and this proves
${\phi_r}|_{\mathcal{C}_X}=
(p_Y|_{\mathcal{C}_Y})^{-1}\circ\overline{r}\circ
p_X|_{\mathcal{C}_X}$.

\subsubsection*{Step 3}
Define a block map
$\Psi_s\colon
\mathcal{B}_{m+1}(\mathcal{C}_Y)\to \mathcal{A}_X$
with $\psi_s$ defined on $\mathcal{C}_Y$ by
$\psi_s(\alpha)_i=\Psi_s
(J_{\alpha(i)}\cdots J_{\alpha(i+m)})$
such that ${\psi_s}\circ \phi_r=Id$
on ${\mathcal{C}_X}$.

We define $\Psi_s$ from $s\colon Y\to X$
just as we defined $\Phi_r$ from $r$ in Step 2.
For $\alpha=\cdots J_{\alpha(-i)} J_{\alpha(-i+1)}
\cdots\in {\mathcal{C}}_Y$
with $J_{\alpha(-i)}\subset \epsilon$,
let $\{J_{\epsilon,l(k)}^{(k)}\}$ be the partition of
$\epsilon \in \mathcal{E}_Y$
for $g^k$, $1\le k \le m+1$, so that
$g^k(J_{\epsilon,l(k)}^{(k)})
=\epsilon_{\alpha(-i+k)} \in \mathcal{E}_Y$
and $g^k(J_{\epsilon,l(k+1)}^{(k+1)})
=J_{\alpha(-i+k)} \subset \epsilon_{\alpha(-i+k)}$.
Then we have a unique descending sequence
\begin{equation}\label{e3}
\epsilon \supset J_{\alpha(-i)}=J_{\epsilon,v(1)}^{(1)}
\supset\cdots \supset J_{\epsilon,v(m+1)}^{(m+1)}
\end{equation}
such that, for $y=p_Y(\alpha)$,
we have $y_i\in J_{\epsilon,v(m+1)}^{(m+1)}$ and 
$s(J_{\epsilon,v(m)}^{(m)})$ is contained in
a unique path $I$ in $X$ such that
$f(I)\in \mathcal{E}_X$.
We define
${\Psi_s}\colon
\mathcal{B}_{m+1}(\mathcal{C}_Y)\to \mathcal{A}_X$
by 
\[
J_{\alpha(-i)}\cdots J_{\alpha(-i+m)} \mapsto I.
\]

Let $p_X|_{\mathcal{C}_X}(a)=x$,
$\overline{r}(x)=y$, and
$(p_Y|_{\mathcal{C}_Y})^{-1}(y)=\alpha$.
If $x_{i+m}\in I_{a(-i-m),u(m+1)}^{(m+1)}$
and $y_i\in J_{\epsilon,v(m+1)}^{(m+1)}$ as in the 
equations (\ref{e2}) and (\ref{e3}),
then we have
$J_{\epsilon,v(m+1)}^{(m+1)}\subset
r(I_{a(-i-m),j(m+1)}^{(m+1)})$
and
\[
s(J_{{\alpha(-i)},v(m+1)}^{(m+1)})\subset
s\circ r(I_{a(-i-m),j(m+1)}^{(m+1)})
=I_{a(-i)}\in \mathcal{A}_X.
\]
Therefore ${\psi_s}\circ \phi_r
\colon \mathcal{C}_X\to \mathcal{C}_X$
is the sliding block code
with memory $m$ and anticipation $m$
induced by a block map defined by
\[
I_{a(-i-m)}\cdots I_{a(-i)}\cdots I_{a(-i+m)}
\mapsto I_{a(-i)},
\]
and ${\psi_s}\circ \phi_r=Id$
on ${\mathcal{C}_X}$.

\subsubsection*{Step 4}
Deduce that $\phi_r$ gives
the required conjugacy $\tilde{r}$.

Because $\phi_r$ maps $\mathcal{C}_X$ onto
$\mathcal{C}_Y$,
and these sets are dense in $\Sigma_X$ and $\Sigma_Y$,
it follows that $\phi_r$ maps 
$\Sigma_X$ onto $\Sigma_Y$.
Similarly $\psi_s$ maps $\Sigma_Y$ onto $\Sigma_X$.
Since the continuous maps $Id|_{\Sigma_X}$ and
$\psi_s\circ \phi_r$ agree on the dense set
$\mathcal{C}_X$, we have
$Id|_{\Sigma_X}=\psi_s\circ \phi_r$ on $\Sigma_X$,
and so $\phi_r$ is a conjugacy.

That $\phi_r$ is a lift of $\overline{r}$ follows
because $p_Y\circ \phi_r=\overline{r}\circ p_X$
on the dense set $\mathcal{C}_X$.
The lifting is unique for it is uniquely determined
on the dense set $\mathcal{C}_X$.
\end{proof}

\begin{remark} \label{3.6}
It is necessary to assume the shift equivalence 
by graph maps.
See Examples \ref{4.3} and \ref{5.a}.
\end{remark}

\section{Graph algorithm}\label{S4}

Suppose that $(X,f)$ is a presentation of a
solenoid satisfying all six Axioms.
Given a finite subset $\mathcal{O}$ of $X$
such that $f(\mathcal{O})=\mathcal{O}$,
we will construct a new presentation
$(X_{\mathcal{O}},f_{\mathcal{O}})$ such that
$(\overline{X}_{\mathcal{O}},\overline{f}_{\mathcal{O}})$
is topologically conjugate to
$(\overline{X},\overline{f})$.
For this purpose,
we will give a graph algorithm
which takes the given presentation $(X,f)$
and $\mathcal{O}$
to produce a new presentation
$(X_{\mathcal{O}},f_{\mathcal{O}})$
and shift equivalence maps
$\rho_{\mathcal{O}}\colon
(X_{\mathcal{O}},f_{\mathcal{O}})\to (X,f)$
and 
$\psi_{\mathcal{O}}\colon
(X,f)\to (X_{\mathcal{O}},f_{\mathcal{O}})$.

\begin{notation}\label{4.n}
By a {\it  path} or {\it directed path} $l$, we  mean
an equivalence class of
locally one-to-one continuous maps
$\gamma\colon [0,1] \to X$
where $\gamma_1\sim \gamma_2$ if and only if
there is an order-preserving homeomorphism
$h_{1,2}\colon [0,1]\to [0,1]$
such that $\gamma_1(t)=\gamma_2\circ h_{1,2}(t)$
for all $t\in [0,1]$.
So a path $l=[\gamma]$ has
the initial point $\gamma(0)$ and
the terminal point $\gamma(1)$.
By abuse of notation,
when we say that a path 
contains a set, we mean that the image of 
the path contains the set.
\end{notation}

Let $\mathcal{P}$ be the set of directed paths 
$l\subset X$ such that
the boundary points of $l$ are
contained in $\mathcal{O}$,
and the interior of $l$ does not contain any point of
$\mathcal{O}$.
Let $l$ be a directed path in $\mathcal{P}$.
Since $f(\mathcal{O})= \mathcal{O}$,
if we consider $f(l)$ as a continuous map
$L\colon [0,1]\to X$,
then the interval $[0,1]$ can be represented
as a union of subintervals
$[0,a_1]\cup [a_1,a_2]\cup\cdots\cup [a_{n-1},1]$
such that $L(a_i)\in \mathcal{O}$ and 
$L\big((a_{i-1},a_{i})\big)\cap\mathcal{O}=\emptyset$
for every $0<i<n$.
Hence for the collection of $[a_{i-1},a_{i}]$
such that
$L([a_{i-1},a_{i}])=l_i\in \mathcal{P}$,
$f\colon X\to X$ induces a {\it wrapping rule}
$\tilde{f}_{\mathcal{P}}\colon
\mathcal{P}\to \mathcal{P}^*$ defined by
$l\mapsto l_{1^\prime}\cdots l_{n^\prime}$.
We write $l_i\in  \tilde{f}(l)$ if
$l_i\in \mathcal{P}$ is
one of these factors of $\tilde{f}(l)$.
(We remark that some factors of $\tilde{f}(l)$ 
may not be paths, because $\tilde{f}(l)$ 
need not be locally one-to-one.)

Our first task is
to find a minimal set of directed paths 
$\mathcal{P}_{\mathcal{O}}$ such that 
$\mathcal{P}_{\mathcal{O}}$ is a finite subset
of $\mathcal{P}$,
$X\subset \bigcup_{l\in\mathcal{P}_{\mathcal{O}}}l$,
$\tilde{f}(l)\in\mathcal{P}_{\mathcal{O}}^*$
for every $l\in \mathcal{P}_{\mathcal{O}}$,
and there exists a positive integer $k$ such that
$l_1 \in \tilde{f}^k(l_2)$
for all $l_1, l_2\in \mathcal{P}_{\mathcal{O}}$.

Let's denote $\mathcal{I}_m$ as 
the set of directed paths
whose boundary points are contained in
$f^{-m}(\mathcal{O})$, and whose interior
does not have any point in $f^{-m}(\mathcal{O})$.
Since we assumed all six Axioms,
there is a positive integer $n$ such that
$f^n(e)\supseteq X$ for every edge $e$ of $X$
from Lemma \ref{2.19}.
So each edge $e$ contains at least one point of 
$f^{-n}(\mathcal{O})$,
and $\mathcal{I}_m$ is a finite set for $m\ge n$.

We have the following lemma
from the Flattening Axiom.

\begin{lemma}\label{4.aaa}
There exists a positive integer $N\ge n$ such that
the interior of each edge contains
at least one point of $f^{-N}(\mathcal{O})$,
and for each vertex $v\in X$
\begin{enumerate}
\item[(1)]
if $v\notin f^{-N}(\mathcal{O})$,
then there exists a path $l_v\in \mathcal{P}$,
unique up to direction, such that
for every path $I \in \mathcal{I}_N$
which contains $v$ as an interior point,
either $f^N(I)\notin \mathcal{P}$
or $f^N(I)=l_v^{\pm 1} \in \mathcal{P}$
where $\pm 1$ denotes the direction, and
\item[(2)]
if $v\in f^{-N}(\mathcal{O})$, then
there exist paths $l_i, l_t\in \mathcal{P}$,
unique up to direction,
such that
if $J_1$ and $J_2$ are elements of $\mathcal{I}_N$
such that $v$ is the terminal point of $J_1$
and the initial point of $J_2$,
then $f^N(J_1J_2)=(l_il_t)^{\pm 1}$. 
\end{enumerate}
\end{lemma}

Fix $N$ satisfying the statement of Lemma \ref{4.aaa}.  
Let
$\mathcal{I}=\{ I\in \mathcal{I}_N \mid
f^N(I)\in \mathcal{P}\}$.
Then $I\in \mathcal{I}_N\backslash \mathcal{I}$
if and only if $f^N(I)$ is not locally one-to-one.
Each $l\in \mathcal{P}$
can be represented as a product of
$I_{i}\in \mathcal{I}_N$ such that
the initial point of $I_{1}$ is
the initial point of $l$,
the terminal point of $I_{i}$ is
the initial point of $I_{i+1}$
for $1\le i < j(l)$, and
the terminal point of $I_{j(l)}$
is the terminal point of $l$
so that
some $I_i\in \mathcal{I}$,
some $I_j\in \mathcal{I}_N\backslash\mathcal{I}$,
and $\tilde{f}^N(l)
=f^N(I_{i(1)})\cdots f^N(I_{i(l)})$
where $I_{i(k)}\in \mathcal{I}$.
Therefore each factor of $\tilde{f}^{N+i}(l)$
which is a path in $X$
is an element of $f^i(\mathcal{I})$ 
for every $l\in\mathcal{P}$ and $i\ge 0$.

\begin{lemma}\label{4.aa}
Suppose $l=f^N(I)\in \mathcal{P}$
for some $I\in \mathcal{I}$.  
Then every factor of $\tilde{f}^i(l)$ is an element of
$f^N(\mathcal{I})$ for every nonnegative integer $i$. 
\end{lemma}
\begin{proof} 
Clearly every factor of $\tilde f^{(i)}(l)$ 
is an element of $f^N(\mathcal I_N)$. 
We must check that every factor is 
locally one to one. 
First consider the case that the image of 
$I$ is a subset of an edge of $X$.
Assume that
$I$ is represented as a product $J_1\cdots J_{j(I)}$
such that each  $J_j\in \mathcal{I}_{N+i}$,
and $f(J_j)\notin \mathcal{I}$ for some $j$. 
Then $f(J_j)\notin \mathcal{I}$ implies that
$f^{N+i}(J_j)$ is not locally one-to-one
on the image of $J_j$,
and so $f^{N+1}(I)$ is not locally one-to-one.
This contradicts the Nonfolding Axiom
as we chose the image of $I$ to be a subset of an edge. 
So we have  $f(J_j)\in \mathcal{I}$
for  $1\le j\le j(I)$.

Now suppose that $I$ contains
a vertex $v$ of $X$ as an interior point.
Let $I^\prime\in \mathcal{I}$ be a subset of an edge.
Then there is a positive integer $k$ 
and a factor $J\in \mathcal{I}_{N+k}$ of $I^\prime$
such that $X\subset f^k(I^\prime)$
and $f^k(J)$ contains $v$ as an interior point.
By Lemma \ref{4.aaa},
we have $f^N(I)=f^{N+k}(J)^{\pm 1}\in \mathcal{P}$,
and $f^{N+i}|_{I}=f^{N+k+i}|_{J}$ is locally
one-to-one as $J$ is a subset of an edge.
Therefore factors of $\tilde{f}^i(f^N(I))$
are elements of $f^N(\mathcal{I})$
for all $i\ge 0$. 
\end{proof}

\begin{definition}[Closed finite path set up to
direction]\label{4.def}
The directed paths which are elements of $\mathcal I$ 
come in pairs, where one path in a pair is 
the other with reversed direction. 
Make a choice of one path from each pair and 
let $\mathcal I_{or}$ be the collection of 
chosen directed paths. 
Define
\[
\mathcal{P}_{\mathcal{O}}
=\{ f^N(I)\mid  I\in \mathcal{I}_{or} \}.
\]
Then $\mathcal{P}_{\mathcal{O}}$ is a finite subset
of $\mathcal{P}$ as $\mathcal{I}$ is a finite set.
\end{definition}

\begin{proposition} \label{4.1}
The set $\mathcal{P}_{\mathcal{O}}$ is
the unique, up to the choice of direction,
minimal subset of $\mathcal{P}$ satisfying  
the following conditions: 
$\tilde{f}(l)\in\mathcal{P}_{\mathcal{O}}^*$
for every $l\in \mathcal{P}_{\mathcal{O}}$,
and $X\subset \bigcup_{l\in\mathcal{P}_{\mathcal{O}}}l$.
There exists a positive integer $k$ such that
$l_1 \in \tilde{f}^k(l)$
for all $l_1, l\in \mathcal{P}_{\mathcal{O}}$.
\end{proposition}
\begin{proof}
By Lemma \ref{4.aa},
$\tilde{f}(l)$ or $\tilde{f}(l)^{-1}$ is contained in
$\mathcal{P}_{\mathcal{O}}^*$
for every $l\in \mathcal{P}_{\mathcal{O}}$.
It remains to check minimality.

Suppose that $l=f^N(I)$ for some $I\in \mathcal{I}$.
Then by Lemma \ref{2.19},
there exists a positive integer $j\ge N$
such that $X\subset f^j(I)$.
If $l_1=f^N(I_1)$ such that
$I_1\in \mathcal{I}$ and
the interior of $I_1$ is contained in an edge $e_1$,
then there exists a subpath $J_1\subset I$ such that
$f^j(J_1)=I_1^{\pm 1}$,
and we have $f^N(I_1)=l_1^{\pm 1}$
is a factor of $\tilde{f}^j(f^{N}(I))$.

Next suppose that $l_2=f^N(I_2)$
where $I_2$ is a path in $\mathcal{I}$ such that
$I_2$ contains a vertex $v$ as an interior point.
Let $I_3\in \mathcal{I}$ be contained in an 
edge $e$,
and $l_3=f^N(I_3)$.
Then for some $m > 0$, $v$ is
the image  under $f^m$ of an interior point of $e$,  
and by Lemma \ref{4.aaa}, 
$l_v^{\pm 1}$ is a factor of $\tilde{f}^{N+m}(l_3)$.
This proves
the minimality of $\mathcal{P}_{\mathcal{O}}$,
and the uniqueness claim is also clear.

It is clear that 
$X\subset \bigcup_{l\in\mathcal{P}_{\mathcal{O}}}l$,
and for all $l_1, l_2\in \mathcal{P}_{\mathcal{O}}$
there exists a positive integer $k=k(1,2)$ such that
$l_1^{\pm 1}$ is a factor of $\tilde{f}^k(l_2)$.
Then the number $k$ can be chosen uniformly
because for every $I\in \mathcal{I}$ such that
the image of 
$I$ is contained in an edge,
if $l=f^N(I)$, then $l^{\pm 1}$ is a factor of
$f^m(l^\prime)$ for every
$l^\prime\in \mathcal{P}_{\mathcal{O}}$ and large $m$.
\end{proof}

\begin{definition}[Construction of new presentation]
\label{4.de}
The new directed graph $X_{\mathcal{O}}$
defined by the set $\mathcal{O}$
has $n$ vertices and $m$ edges
where $n$ is the cardinality of $\mathcal{O}$ and
$m$ is the cardinality of $\mathcal{P}_{\mathcal{O}}$. 
The set of vertices $\mathcal{V}_{\mathcal{O}}$
of $X_{\mathcal{O}}$ corresponds to $\mathcal{O}$,
and the set of edges $\mathcal{E}_{\mathcal{O}}$
corresponds to $\mathcal{P}_{\mathcal{O}}$
by the following rule:
If $l_i \in \mathcal{P}_{\mathcal{O}}$ is
a directed path from $v_{i,1}$ to $v_{i,2}$
in $X$ represented by a continuous map
$\gamma_{l_i}\colon [0,1]\to X$
as in Notation \ref{4.n},
then $e_i\in \mathcal{E}_{\mathcal{O}}$
is a directed edge from $v_{i,1}$ to $v_{i,2}$
represented by a continuous map
$\gamma_{i}\colon [0,1]\to e_i$
such that 
$\gamma_i|_{(0,1)}$ 
is a homeomorphism,
and $l_i^{-1}$ corresponds to $e_i^{-1}$
defined by $\gamma_{i}^{-1}\colon [0,1]\to e_i$
such that $\gamma_{i}^{-1}(t)=\gamma_i(1-t)$.
So there is a natural projection
$\rho_{\mathcal{O}}\colon X_{\mathcal{O}}\to X$
defined by
$\gamma_{l_i}\circ (\gamma_{i})^{-1}\colon
e_i \to [0,1] \to l_i$,
and the graph map $f_{\mathcal{O}}\colon
X_{\mathcal{O}}\to X_{\mathcal{O}}$
is induced by $f\colon X\to X$ satisfying
$\rho_{\mathcal{O}}\circ f_{\mathcal{O}}=
f\circ \rho_{\mathcal{O}}$.
Hence if
$\tilde{f}_{\mathcal{P}}\colon
l_i\mapsto l_{i,1}^{s(1)}\cdots l_{i,m(i)}^{s(m(i))}
\in \mathcal{P}_{\mathcal{O}}^*$
for $l_i\in \mathcal{P}_{\mathcal{O}}$, then
$\tilde{f}_{\mathcal{O}}\colon e_i\mapsto
e_{i,1}^{s(1)}\cdots e_{i,m(i)}^{s(m(i))}
\in \mathcal{E}_{\mathcal{O}}^*$
where $\rho_{\mathcal{O}}(e_i)=l_i
\in \mathcal{E}_{\mathcal{O}}$
and $s(i)$ denotes the direction.
\end{definition}

\begin{remark}\label{4.rr}
Suppose that
$\mathcal{P}_{\mathcal{O}}^\prime$ is
$\mathcal{P}_{\mathcal{O}}$ with
a different choice of directions of paths.
Let $(X_{\mathcal{O}},f_{\mathcal{O}})$
be defined by $\mathcal{P}_{\mathcal{O}}$, and
$(X_{\mathcal{O}}^\prime,f_{\mathcal{O}}^\prime)$
 defined by $\mathcal{P}_{\mathcal{O}}^\prime$.
Then the two graphs are the same except that 
the directions of some corresponding edges 
might be reversed, and it  
 is easy to see that 
$f_{\mathcal{O}}$
and $f_{\mathcal{O}}^\prime$
are shift equivalent by graph maps. 
\end{remark}

\begin{example} \label{4.3}
Let $X$ be a wedge of two circles $a,b$
with a unique vertex $p$, and
$f\colon X\to X$ defined by
$a\mapsto aab$ and $b\mapsto ab$.
So $(X,f)$ is given by the following diagram,
in which $p$ is the vertex of $X$
and $\{q,r\}$ is a period $2$ orbit.
\begin{figure}[ht]
\centerline{\scalebox{.195}{\includegraphics{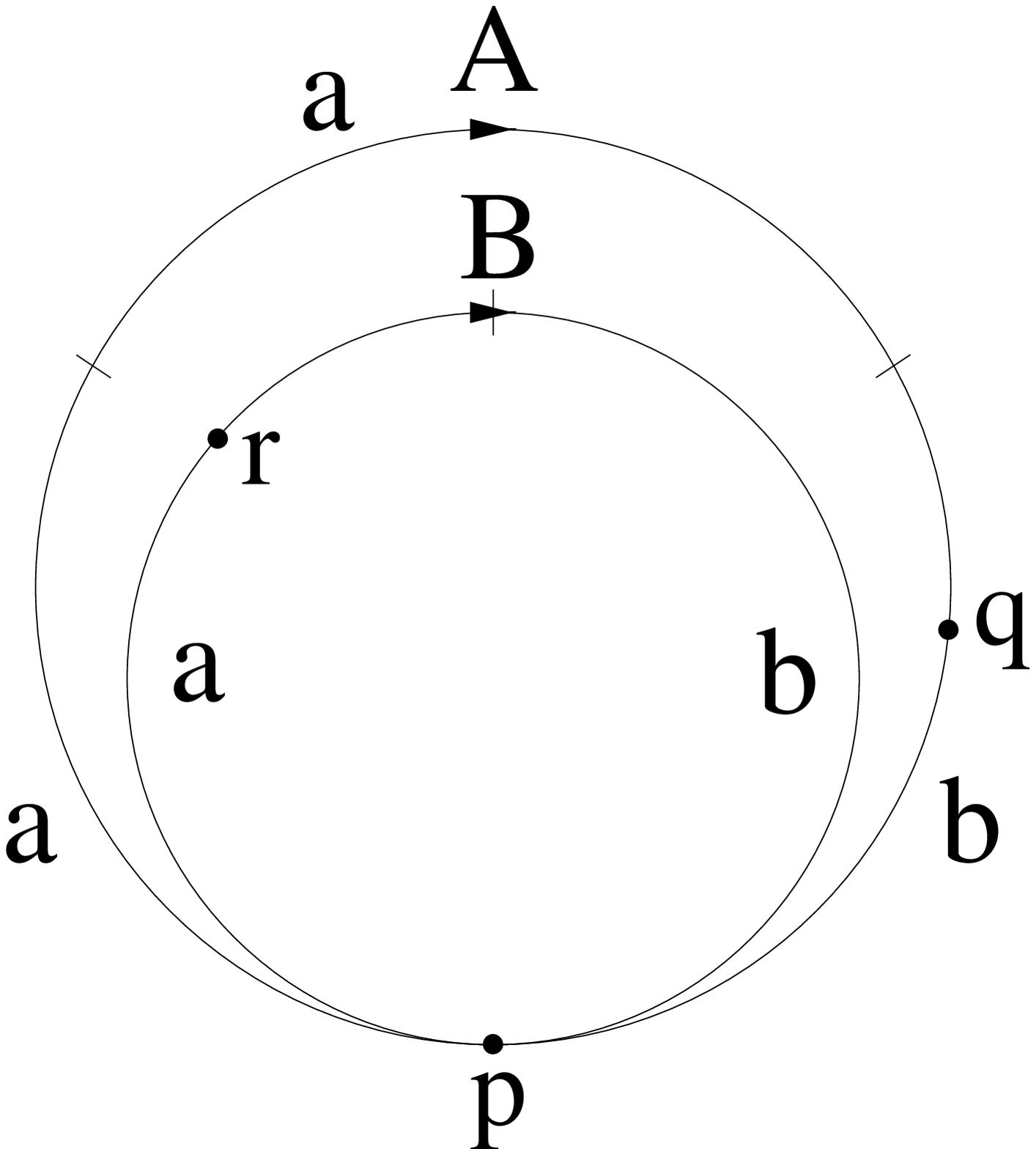}}}
\end{figure}

The set of directed paths $\mathcal{P}$
defined by the indicated periodic orbit $\{q,r \}$
is
$\{ \alpha,\beta,\gamma,\delta, \epsilon,\zeta,
\alpha^{-1},\beta^{-1},\gamma^{-1},\delta^{-1},
\epsilon^{-1},\zeta^{-1}\}$
where $\alpha$ is the circle $A$ based at $q$,
$\beta$ is the path from $q$ through $p$ to $r$,
$\gamma$ is from $r$ through $p$ to $q$,
$\delta$ is the circle $B$ based at $r$,
$\epsilon$ is the path from $q$ through $p$
to $r$ with reverse direction to $B$, and
$\zeta$ is the path from $r$ through $p$
with reverse direction to $B$ to $q$.

Every edge of $X$ has at least two points
of $f^{-1}({\{q,r}\})$,
and it is not difficult to check that
$f(\mathcal{I})=\{\alpha, \beta, \gamma,
\alpha^{-1}, \beta^{-1}, \gamma^{-1}\}$.
So $\mathcal{P}_{{\{q,r}\}}$ is 
$\{\alpha,\beta,\gamma\}$
up to the choice of direction,
and the induced wrapping rule
$\tilde{f}_{\mathcal{P}}\colon
         \mathcal{P}\to \mathcal{P}^*$
is given by
\[
\alpha\mapsto \gamma\alpha\beta,\quad
\beta\mapsto \gamma,\quad
\gamma\mapsto \beta\gamma\alpha\beta.
\]
Hence 
the new graph $X_{\{q,r\}}$ defined by $\{q,r\}$ is
the following graph.
\begin{figure}[!h]
\centerline{\scalebox{.35}{\includegraphics{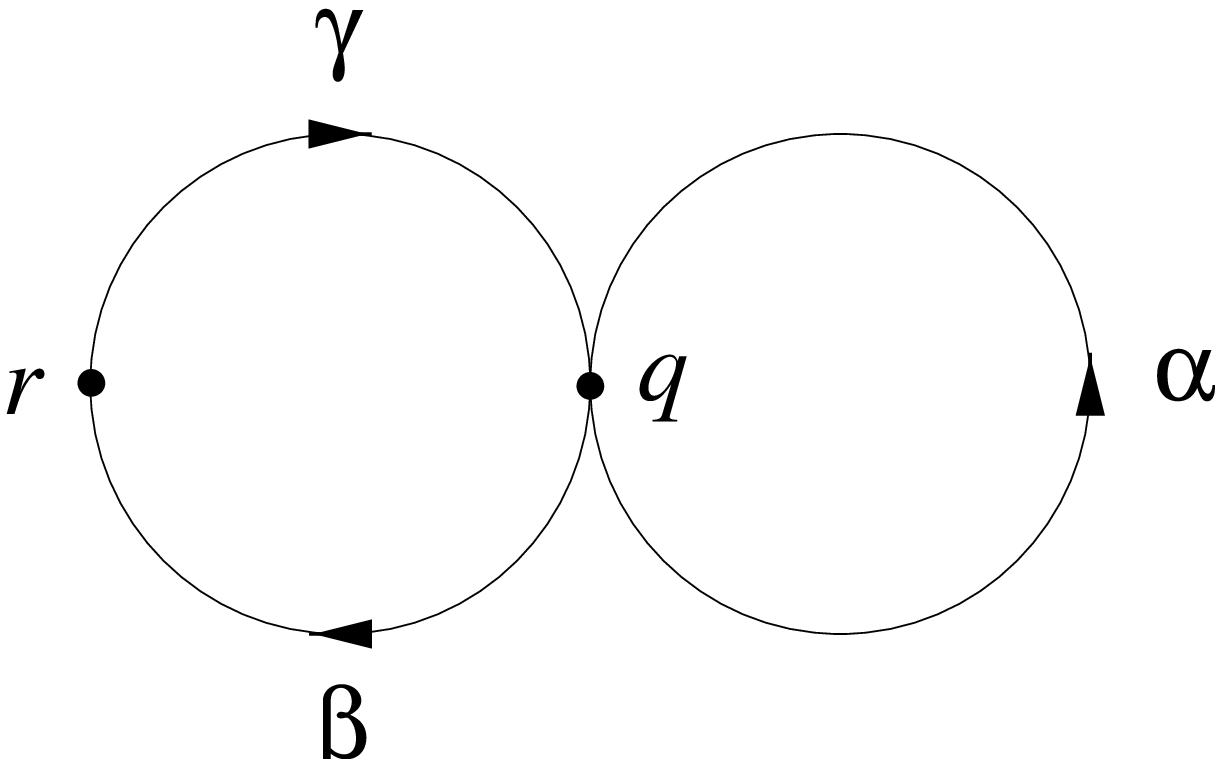}}}
\end{figure}
\end{example}

\begin{remark}
To compute the graph $X_{\mathcal O}$ 
and the map $f_{\mathcal O}$, 
we don't need to find
the integer $N$ and $\mathcal{I}$
given in Lemma \ref{4.aaa}.
If we choose a path $l\in \mathcal{P}$,
and iteratively
 apply $\tilde{f}$ to the factors of
$\tilde{f}^n(l)$ which are paths, 
we will eventually generate 
a set of paths invariant under 
$\tilde f$, and this set will contain 
$\mathcal P_{\mathcal O}$. 
\end{remark}

\begin{example} \label{4.4}
Suppose that $(Y,g)$ is given by Figure \ref{xcv6},
and $p$ is a fixed point of $g$.
\begin{figure}
\centerline{\scalebox{.27}{\includegraphics{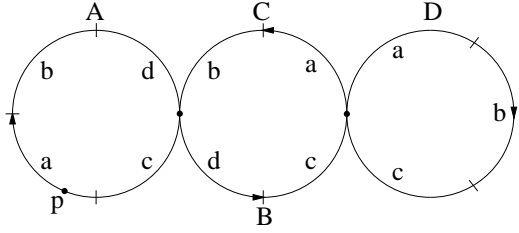}}}
\caption{$(Y,g)$ with a fixed point $p$}\label{xcv6}
\end{figure}
Then the set $\mathcal{P}$ of directed paths
based at $p$ is an infinite set
for if we call $a_1$ the path
from $p$ to the branch point in the edge $a$
and $a_2$ the path from the branch point to $p$, then
the paths $a_1b\underbrace{d\cdots d}_{n}ca_2$
are in $\mathcal{P}$.

If we choose $\ell$ as the loop $a_1a_2$ based at $p$,
then $\tilde{g}(\ell)=a_1bdca_2$ and
$\tilde{g}^2(\ell)=a_1bddca_2\, a_1bca_2\, a_1bca_2\,
a_1bca_2$.

For $e=a_1bddca_2$ and $f=a_1bca_2$,
$\tilde{g}\colon e\mapsto efff
\text{ and } f\mapsto ef$.
So $\{e,f\}$ is a closed minimal subset of
$\mathcal{P}$, and
$\mathcal{P}_{\{p\}}=\{e,f\}$ by uniqueness.
The new graph $(Y_{\{p\}},g_{\{p\}})$
defined by $\{p\}$ is the following graph.
\begin{figure}[h]
\centerline{\scalebox{.19}{\includegraphics{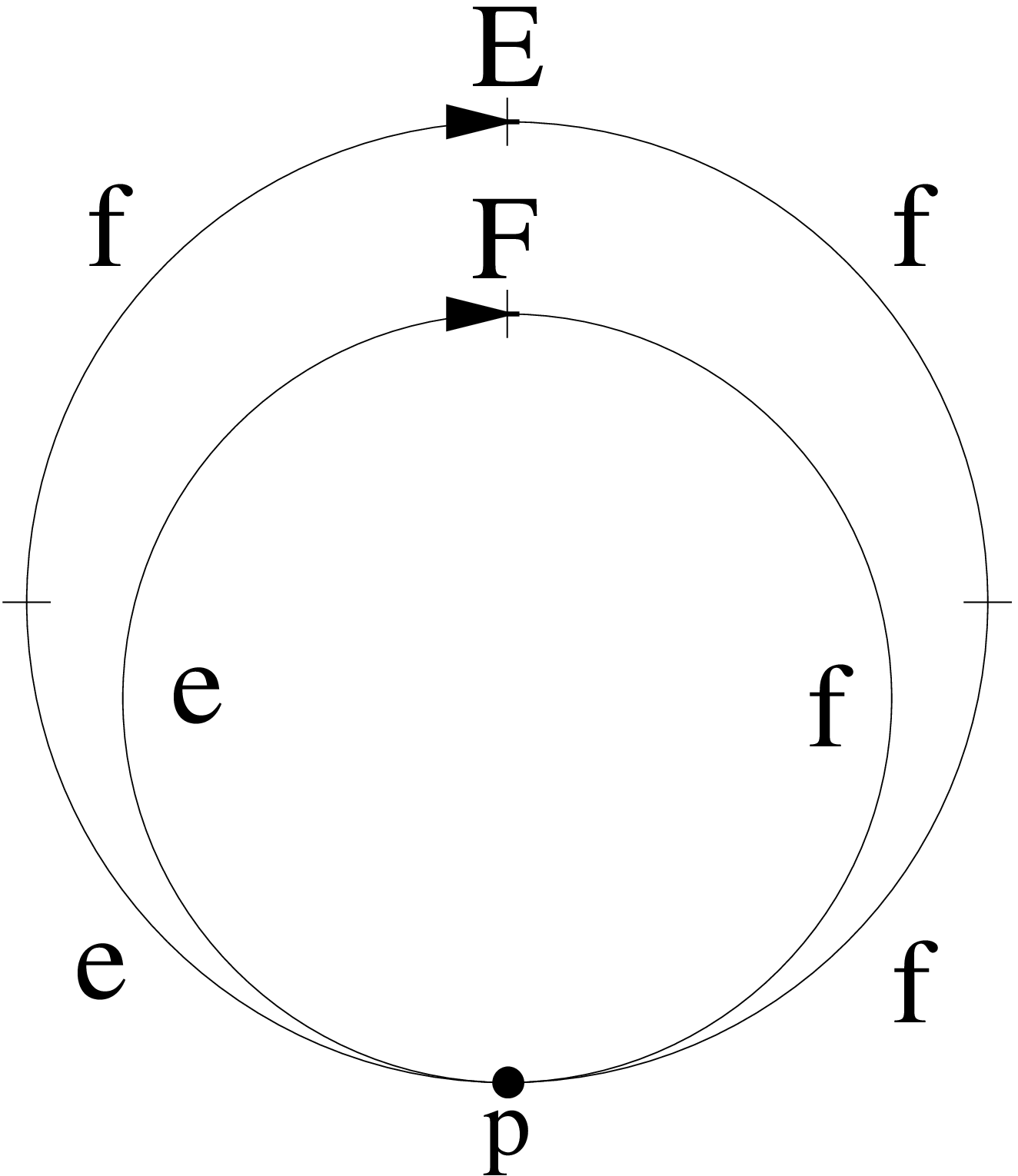}}}
\end{figure}
\end{example}

\begin{example}\label{4.e}
Suppose that $Z$ is given in the following graph,
that $h\colon Z \to Z$ is given by
$a\mapsto bab^{-1}$ and $b \mapsto aba^{-1}$,
and that $p$ is a fixed point.
\begin{figure}[!h]
\centering \includegraphics[height=0.8cm]{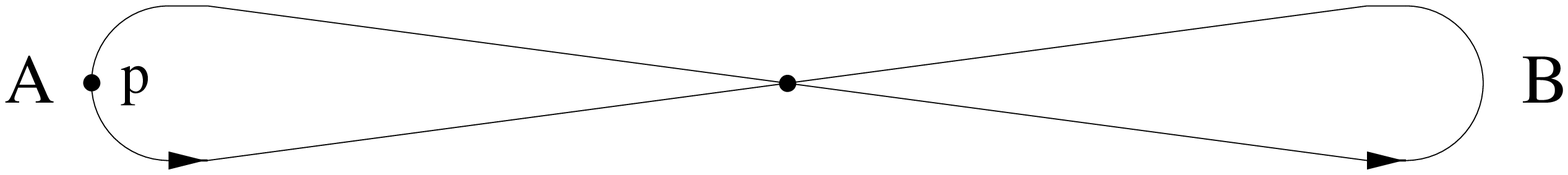}
\end{figure}

Let $a_1$ denote
the path from $p$ to the branch point,
$a_2$ the path from the branch point to $p$,
and $\ell=a_2ba_1$.
Then
$\tilde{h}(\ell)=a_2b^{-1}a_1\, a_2ba_2^{-1}\,
a_1^{-1}ba_1$ and for
$\alpha=a_2b^{-1}a_1$, $\beta=a_2ba_2^{-1}$
and $\gamma=a_1^{-1}ba_1$,
\[
\tilde{h}\colon
\alpha\mapsto \alpha\beta^{-1}\gamma,\quad
\beta\mapsto \alpha\beta\alpha^{-1},\quad
\gamma\mapsto \gamma^{-1}\beta\gamma.
\]
The new graph $Z_{\{p\}}$ defined by $\{p\}$
is a wedge of three circles as in the graph below. 
\begin{figure}[h]
\centering \includegraphics[height=1.8cm]{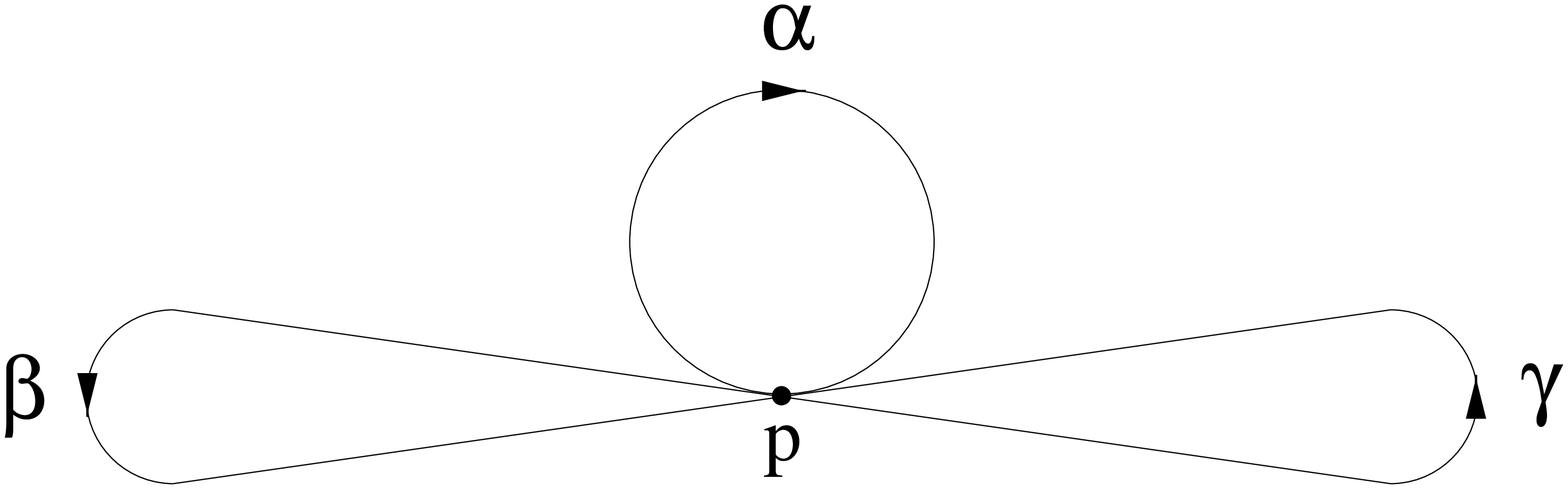}
\end{figure}
\end{example}

\begin{theorem} \label{4.5}
Suppose that $(X,f)$ satisfies all six Axioms,
that $\mathcal{O}$ is a finite subset of $X$
such that $f(\mathcal{O})=\mathcal{O}$,
and that $\mathcal{P}_{\mathcal{O}}$ and
$(X_{\mathcal{O}},f_{\mathcal{O}})$
are the  minimal closed subset of $\mathcal{P}$
and the new presentation defined by $\mathcal{O}$,
respectively.
Then the natural projection
$\rho_{\mathcal{O}}\colon X_{\mathcal{O}}\to X$
gives a conjugacy
from 
$(\overline{X}_{\mathcal{O}},\overline{f}_{\mathcal{O}})$
to $(\overline{X},\overline{f})$.
\end{theorem}
\begin{proof}
We will show that $\rho=\rho_{\mathcal{O}}$
is part of a shift equivalence, that is,
we will define a continuous map
$\psi\colon X \to X_{\mathcal{O}}$ 
and a positive integer $m$ such that
\[
f_{\mathcal{O}}\circ \psi=\psi\circ f,\text{ }
\rho \circ \psi = f^m, \text{ and }
\psi \circ \rho =f_{\mathcal{O}}^m.
\]

Recall 
 $\mathcal{P}_{\mathcal{O}} =  
\{f^N(I)\mid I\in \mathcal I_{or} \}$
for some fixed choice of a subset $\mathcal I_{or}$ 
of $\mathcal I_N$ as in Definition \ref{4.def}.
By the choice of $N$ (from Lemma \ref{4.aaa}),  
every $I\in \mathcal{I}_{or}$
contains at most one vertex of $X$.
It is trivial that
if $x\in X\backslash f^{-N}(\mathcal{O})$
is contained in some $I\in \mathcal{I}_{or}$ and 
the interior of $I$ does not contain any vertex of $X$,
then $I$ is the unique element in 
$\mathcal{I}_{or}$ which contains $x$.

If $v$ is a vertex of $X$, then by Lemma \ref{4.aaa},
$v\in f^{-N}(\mathcal{O})$ or
there is a neighborhood $U_v$ of $v$ such that
$U_v$ is the union of all paths in $\mathcal{I}$
containing $v$ as an interior point and
$f^{N}(U_v)$ is the image of
$l_v\in \mathcal{P}_{\mathcal{O}}$.
Hence,
for every $x\in X\backslash f^{-N}(\mathcal{O})$,
there is a unique path
$l_x\in \mathcal{P}_{\mathcal{O}}$
such that 
if $x\in I\in \mathcal I_{or}$, 
then $f^N(I) = l_x^{\pm 1}$.

For $x \in X\backslash f^{-N}(\mathcal O)$, let 
$e_x$ be the edge of $X_{\mathcal{O}}$
corresponding to $l_x$.
We will define $\psi (x)$ to be the appropriate 
point $x_{\mathcal O}$ in $e_x$ 
satisfying $\rho\circ \psi (x) = f^N(x)$. 
Let $\gamma \colon [0,1] \to e_x$ be the 
continuous function (in the equivalence class $l_x$) 
associated to $e_x$ 
in the definition of $X_{\mathcal O}$
and  $\rho$. 
Fix $I$ in $\mathcal I_{or}$ such that $x\in I$.
For the moment
let $I$ also denote a specific map $[0,1] \to X$. 
Then there is a homeomorphism
$h\colon [0,1]\to [0,1]$ such that 
$\rho \circ \gamma = f^N \circ I \circ h$. 
Let $t$ be the unique number in $(0,1)$
such that $I\circ h(t) = x$,
and define $x_{\mathcal O} = \gamma (t)$. 
Then $\rho\circ\psi(x)=\rho(x_{\mathcal O})  
=\rho\circ\gamma(t)=f^N\circ I\circ h(t)=f^N(x)$ 
as required. 

For $x\in f^{-N}(\mathcal O)$,
we define $\psi(x)$ as the unique point
in $X_{\mathcal O}$ which $\rho$ maps to
$f^N(x)\in \mathcal O$. 
Then $\psi$ is continuous, $\rho \circ \psi = f^N$, and
clearly $f \circ \rho=\rho\circ f_{\mathcal{O}}$
by Definition \ref{4.de}.
Consequently
$\rho\circ \psi\circ \rho = f^N\circ\rho
=\rho \circ f_{\mathcal O}^N$. 
Therefore the two maps
$\psi\circ \rho$ and $f_{\mathcal O}^N$ send
any given edge (considered as a path) in $X_{\mathcal O}$ 
to paths which $\rho$ sends to the same  
concatenation of elements of $\mathcal I$.  
Such a concatenation has a unique lifting under $\rho$, 
therefore $\psi\circ \rho = f_{\mathcal O}^N$.

It remains to show that
$f_{\mathcal{O}}\circ \psi=\psi\circ f$.
Because $\rho$ is surjective, it suffices to show
$f_{\mathcal{O}}\circ \psi\circ\rho
=\psi\circ f\circ\rho$, 
and this is true because 
\[
f_{\mathcal{O}}\circ \psi\circ \rho
=f_{\mathcal{O}}^{N+1}
= (\psi\circ \rho )\circ f_{\mathcal{O}}
=\psi \circ (\rho \circ f_{\mathcal{O}})
=\psi\circ (f \circ \rho).
\] 
Therefore $(X_{\mathcal{O}},f_{\mathcal{O}})$ is 
shift equivalent to $(X,f)$, and
$\overline{\rho}_{\mathcal{O}}$ is
a topologically conjugacy by Proposition \ref{3.4}.
\end{proof}

\begin{remarks} \label{4.r}
\begin{enumerate}
\item[(1)]
Theorem \ref{4.5} requires the Flattening Axiom.
See Example \ref{4.10}.
\item[(2)]
If $\mathcal{O}$ is not a subset of
$\mathcal{V}_X$,
then the shift equivalence maps $\rho$ and $\psi$
in Theorem \ref{4.5} cannot be graph maps.  
\end{enumerate}
\end{remarks}

\begin{remark}[Preperiodic sets]
If $\mathcal{F}$ is a finite subset of $X$
such that $f(\mathcal{F})\subseteq \mathcal{F}$,
then it is not difficult to apply
the graph algorithm to $\mathcal{F}$
so that we have a finite graph with a graph map
$(X_{\mathcal{F}},f_{\mathcal{F}})$.
And there is a positive integer $k=k_{\mathcal{F}}$
such that
$f^k(\mathcal{F})=f^{k+i}(\mathcal{F})$
for all $i\ge 0$. 
Let $\mathcal{O}_{\mathcal{F}}=f^k(\mathcal{F})$.

We get $f$ and $f_{\mathcal{F}}$ shift equivalent
just as before.
Now it is not hard to check that
$f_{\mathcal{F}}$ and $f_{\mathcal{O}_{\mathcal{F}}}$
are shift equivalent by graph maps.
Hence the associated SFT covers of the next section
will be conjugate for $\mathcal{F}$ and
$\mathcal{O}_{\mathcal{F}}$,
and that is why we only concern ourselves
with sets $\mathcal{O}$ which are unions of
periodic orbits.
\end{remark}

\subsection*{Elementary presentation}
One interesting application of the graph algorithm
is the elementary presentations of solenoids.
In \cite[\S5]{w2}, Williams introduced an 
{\it elementary presentation} of a solenoid 
in which $X$ is a wedge of circles and 
$f$ leaves the unique branch point of $X$ fixed.
And he showed in \cite[Theorem 5.2]{w2} that, 
for every generalized 1-dimensional solenoid
$(\overline{X}, \overline{f})$, 
there exists  an integer $m$ such that
$(\overline{X},\overline{f^m})$
has an elementary presentation.
We can improve his theorem by getting
sharp bounds on $m$.

Suppose that $(X,f)$ satisfies all six axioms,
and that $a$ is a fixed point of $f^m$.
For $(X,f^m)$, let $(X_{\{a\}},f^m_{\{a\}})$
be the new presentation defined by $\{a\}$.
Then the new graph $X_{\{a\}}$
has only one vertex $a$ 
which is a fixed point by
$f^m_{\{a\}}\colon X_{\{a\}}\to X_{\{a\}}$,
and each edge in $X_{\{a\}}$ is homeomorphic to a circle.
So $(X_{\{a\}},f^m_{\{a\}})$ is
an elementary presentation
and $(\overline{X}, \overline{f^m})$ is conjugate to
$(\overline{X}_{\{a\}},\overline{f^m}_{\{a\}})$
by Theorem \ref{4.5}.
More generally we have the following proposition:

\begin{proposition} \label{4.8}
For a given 1-solenoid,
the minimal number of vertices 
in a presentation $(X,f)$
is the minimal period of points in $X$.
In particular, $(\overline{X}, \overline{f})$
has an elementary presentation if and only if
$f\colon X\to X$ has a fixed point.
\end{proposition}

\begin{remark}
Williams showed that
two elementary presented solenoids
$(\overline{Y}_1,\overline{g}_1)$
and $(\overline{Y}_2,\overline{g}_2)$
are topologically conjugate to each other
if and only if the shift equivalence classes of
${g_1}_*\colon \pi_1(Y_1, y_1) \to \pi_1(Y_1, y_1)$
and 
${g_2}_*\colon \pi_1(Y_2,y_2) \to \pi_1(Y_2,y_2)$
are the same where $y_i$ is the unique branch point
of $Y_i$ for $i=1,2$ (\cite[7.3]{w2}).
Proposition \ref{4.8} extends the range of
Williams' classification theorem.
\end{remark}

\begin{example} \label{4.10}
Theorem \ref{4.5} and Proposition \ref{4.8}
require the Flattening Axiom.

Let $X$ be a wedge of two circles $a$ and $b$
with $f\colon X\to X$ defined by
\[
a\mapsto bba,\quad b\mapsto abb.
\]
Then $(X,f)$ is an elementary presentation of
a branched solenoid and $f$ does not satisfy 
the Flattening Axiom.
The circle $b$ contains
a unique nonbranch fixed point $q$.

The directed path set $\mathcal{P}_{\{q\}}$
has three loops
\[
\alpha=b_2b_1,\quad \beta=b_2ab_2,\quad \gamma=b_2aab_1
\]
where $b_1$ is the arc from the branch point to $q$
and $b_2$ is from $q$ to the branch point.

Let $X_{\{q\}}$ be a wedge of three circles
$\alpha, \beta, \gamma$ based at $q$,
and $f_{\{q\}}\colon X_{\{q\}}\to X_{\{q\}}$
the map induced from $f\colon X\to X$ by
\[
\alpha\mapsto \alpha\beta,\quad
\beta\mapsto \alpha\alpha\alpha\gamma,\quad
\gamma\mapsto \alpha\alpha\alpha\beta\alpha\gamma.
\]
Then $(X_{\{q\}},f_{\{q\}})$
satisfies the Flattening Axiom.
So $(\overline{X},\overline{f})$ is not
topologically conjugate to
$(\overline{X}_{\{q\}},\overline{f}_{\{q\}})$.
\end{example}

\section{Canonical SFT covers}\label{S5}

In \S4, we constructed a new presentation
$(X_{\mathcal{O}},f_{\mathcal{O}})$
for the given presentation $(X,f)$ and
the finite invariant set $\mathcal{O}$ of $X$.
Associated to this new presentation,
we have an SFT
$(\Sigma_{\mathcal{O}},\sigma_{\mathcal{O}})$
defined by
$(\mathcal{E}_{X_{\mathcal{O}}},\tilde{f}_{\mathcal{O}})$
and a semiconjugacy
$p_{\mathcal{O}}\colon
(\Sigma_{\mathcal{O}},\sigma_{\mathcal{O}})\to
(X_{\mathcal{O}},f_{\mathcal{O}})$
constructed by the standard algorithm described in
section 2.
We will show that this new SFT cover defined by
a finite invariant set is canonical.
As an application, we show that
two solenoids are not conjugate by comparing
Bowen-Franks groups of their SFT covers
defined by periodic orbits of the same period.

Suppose that $(X,f)$ satisfies all six Axioms,
that $\mathcal{O}$ is a finite subset of $X$
such that $f(\mathcal{O})=\mathcal{O}$,
and that $(X_\mathcal{O},f_\mathcal{O})$
is the presentation defined by
$(X,f)$ and $\mathcal{O}$.
If $\mathcal{E}_\mathcal{O}$ is the set of edges
in $X_\mathcal{O}$, and 
$\tilde{f}_\mathcal{O}\colon 
\mathcal{E}_\mathcal{O}\to \mathcal{E}_\mathcal{O}^*$
is the wrapping rule associated to $f_\mathcal{O}$,
then $(\mathcal{E}_\mathcal{O},\tilde{f}_\mathcal{O})$
defines a two-sided SFT 
$(\Sigma_{\mathcal{O}},\sigma_{\mathcal{O}})$.

There is a well-defined quotient map
$p_\mathcal{O}\colon
\Sigma_{\mathcal{O}}\to \overline{X}_{\mathcal{O}}$
such that
$\overline{f}_\mathcal{O}\circ p_\mathcal{O}
=p_\mathcal{O}\circ \sigma_\mathcal{O}$.
If $\rho_\mathcal{O}\colon X_\mathcal{O}\to X$
is the natural projection which maps
each edge $e$ of $X_\mathcal{O}$
to the corresponding path $l$ in $X$,
then
$\overline{\rho}_\mathcal{O}\colon
\overline{X}_\mathcal{O}\to \overline{X}$
is a conjugacy by Theorem~\ref{4.5},
and $\overline{\rho}_\mathcal{O}\circ p_\mathcal{O}
\colon \Sigma_\mathcal{O}\to \overline{X}$
is a finite-to-one quotient map.

The canonical projection map to the
zeroth coordinate 
$\pi\colon \overline{X}\to X$ induces a bijection
$\overline{\mathcal{O}}\leftrightarrow\mathcal{O}$
of finite invariant sets of $\overline{X}$ and $X$.
We will call
$(\Sigma_{\mathcal{O}},\sigma_{\mathcal{O}})$
the {\bf SFT cover of $(\overline{X},\overline{f})$
defined by $\mathcal{O}$} or by $\overline{\mathcal{O}}$.

\begin{example} \label{5.a}
Let $(X,f)$ be as in Example \ref{4.3},
$p$ the branch point of $X$ which is a fixed
point of $f$, and $\{q,r\}$ a period $2$ orbit.
Then the `natural' SFT cover of
$(\overline{X},\overline{f})$
is the SFT $(\Sigma_{\{p\}},\sigma_{\{p\}})$
defined by the orbit $\{p\}$. 
From the wrapping rule
$\tilde{f}\colon a\mapsto aab, b\mapsto ab$
we see that
$(\Sigma_{\{p\}},\sigma_{\{p\}})$ is
represented by the following adjacency matrix.
\[
M_{\{p\}}=
\begin{pmatrix}
2&1\\ 
1&1
\end{pmatrix}
\]

The induced map
$\tilde{f}_{{\{q,r}\}}\colon
\mathcal{P}_{{\{q,r}\}}\to \mathcal{P}_{{\{q,r}\}}^*$
is given by
\[
\alpha\mapsto \gamma\alpha\beta,\quad
\beta\mapsto \gamma,\quad
\gamma\mapsto \beta\gamma\alpha\beta.
\]
So the SFT cover $(\Sigma_{\{q,r\}},\sigma_{\{q,r\}})$
of $(\overline{X},\overline{f})$ defined by
$\{q,r\}$ is given by the following matrix.
\[
M_{\{q,r\}}=
\begin{pmatrix}
1&1&1\\ 
0&0&1\\
1&2&1
\end{pmatrix}
\]
Remark that $(\overline{X}_{\{p\}},\overline{f}_{\{p\}})$
is topologically conjugate to 
$(\overline{X}_{\{q,r\}},\overline{f}_{\{q,r\}})$
by Theorem \ref{4.5}.
But $(\Sigma_{\{p\}},\sigma_{\{p\}})$
is not topologically conjugate to 
$(\Sigma_{\{q,r\}},\sigma_{\{q,r\}})$
as the trace of $M_{\{p\}}$
is different from that of $M_{\{q,r\}}$.
\end{example}

\begin{theorem} \label{5.aa}
Suppose that
$(X,f)$ and $(Y,g)$ satisfy all six Axioms,
and that
$(\overline{X},\overline{f})$ is topologically
conjugate to $(\overline{Y},\overline{g})$
by a conjugacy map $\phi$.
If $\overline{\mathcal{O}}$
is a finite union of periodic orbits of  
$\overline{f}$ and
$\overline{\mathcal{O}^\prime}=
\phi(\overline{\mathcal{O}})$,
then there is a unique conjugacy
$\tilde{\phi}\colon
(\Sigma_{\mathcal{O}},\sigma_{\mathcal{O}}) \to
(\Sigma_{\mathcal{O}^\prime},
         \sigma_{\mathcal{O}^\prime})$
such that 
the following diagram commutes.
\[
\begin{CD}
(\Sigma_{\mathcal{O}},\sigma_{\mathcal{O}})
@>\tilde{\phi}>>
(\Sigma_{\mathcal{O}^\prime},
              \sigma_{\mathcal{O}^\prime})\\
@Vp_{\mathcal{O}}VV       @VVp_{\mathcal{O}^\prime}V\\
(\overline{X}_{\mathcal{O}},\overline{f}_{\mathcal{O}})
@.
(\overline{Y}_{\mathcal{O}^\prime},
            \overline{g}_{\mathcal{O}^\prime})\\
@V\overline{\rho}_{\mathcal{O}}VV 
@VV\overline{\rho}_{\mathcal{O}^\prime}V\\
(\overline{X},\overline{f})
@>>\phi>
(\overline{Y},\overline{g})
\end{CD}
\]
\end{theorem}
\begin{proof}
That $\phi$ is a conjugacy  implies that
$\overline{\mathcal{O}^\prime}$ is a finite union
of periodic orbits of $\overline{g}$.
Let $(X_{\mathcal{O}},f_{\mathcal{O}})$ and 
$(Y_{\mathcal{O}^\prime},g_{\mathcal{O}^\prime})$ 
be the new graphs with graph maps defined by
$\mathcal{O}=\pi(\overline{\mathcal{O}})$ and
$\mathcal{O}^\prime=\pi(\overline{\mathcal{O}^\prime})$,
respectively.
Then by Theorem \ref{4.5},
there exist shift equivalence maps
$(\rho_{\mathcal{O}}, \psi_{\mathcal{O}})$
for $(X_{\mathcal{O}},f_{\mathcal{O}})$ and 
$(X,f)$,
and
$(\rho_{\mathcal{O}^\prime},
\psi_{\mathcal{O}^\prime})$ for
$(Y_{\mathcal{O}^\prime},g_{\mathcal{O}^\prime})$ and 
$(Y,g)$.
Let $\phi_{\mathcal{O}}\colon
\overline{X}_{\mathcal{O}} \to
\overline{Y}_{\mathcal{O}^\prime}$
be the conjugacy 
$(\overline{\rho}_{\mathcal{O}^\prime})^{-1}
\circ \phi \circ \overline{\rho}_{\mathcal{O}}$,
which lifts $\phi$.

By the Ladder Lemma,
there is a shift equivalence 
$r_{\mathcal{O}}\colon
X_{\mathcal{O}} \to Y_{\mathcal{O}^\prime}$
and 
$s_{\mathcal{O}}\colon
Y_{\mathcal{O}^\prime} \to X_{\mathcal{O}}$
such that
$\overline{r}_{\mathcal{O}}=\phi_{\mathcal{O}}$ and
$\overline{s}_{\mathcal{O}}=\phi_{\mathcal{O}}^{-1}$.
Since $\phi$ sends $\overline{\mathcal{O}}$
to $\overline{\mathcal{O}^\prime}$,
$r_{\mathcal{O}}$ and $s_{\mathcal{O}}$
are graph maps.
Then it follows from Theorem \ref{3.5} that
there is a unique conjugacy
$\tilde{\phi}\colon
(\Sigma_{\mathcal{O}},\sigma_{\mathcal{O}})
\to
(\Sigma_{\mathcal{O}^\prime},
   \sigma_{\mathcal{O}^\prime})$
lifting $\phi_{\mathcal{O}}$.
Therefore $\tilde{\phi}$ is the unique conjugacy
lifting $\phi$.
\end{proof}

\begin{remarks} \label{5.b}
\begin{enumerate}
\item[(1)]
It is necessary to assume
$\overline{\mathcal{O}^\prime}
=\phi(\overline{\mathcal{O}})$. See Example \ref{5.a}.
\item[(2)]
We need the Flattening Axiom to guarantee
that $\phi_{\mathcal{O}}$ is a conjugacy.
See Remark \ref{4.r}.
\end{enumerate}
\end{remarks}

\subsection*{Bowen-Franks groups}
We say that two dynamical systems
$(M,\phi)$ and $(M^\prime,\phi^\prime)$
are {\bf flow equivalent}
if they have topologically equivalent suspension flows.

\begin{definition} \label{5.4}
(\cite[\S7]{lm})
Let $A$ be an $r\times r$ integral matrix.
The {\bf Bowen-Franks group} of $A$ is
\[
BF(A)=\text{coker}\,(Id - A)
={\mathbb{Z}^r}/{\mathbb{Z}^r(Id - A)},
\] 
where $\mathbb{Z}^r(Id - A)$ is the image of
$\mathbb{Z}^r$
under the matrix $Id - A$ acting on the right.
\end{definition}

If two irreducible SFTs are flow equivalent,
then they have isomorphic Bowen-Franks group (\cite{bf}).

\begin{example} \label{5.7}
Williams posed the following question
(\cite{us, w2}):
If $X$ is a wedge of two circles $a, b$,
 and $g_1, g_2\colon X\to X$ are given by
\begin{alignat*}{2}
\tilde{g}_1(a)&=aabba,&\quad \tilde{g}_1(b)&=a,
\text{ and}\\
\tilde{g}_2(a)&=ababa,  &\quad \tilde{g}_2(b)&=a,
\end{alignat*}
then are $g_1$ and $g_2$ shift equivalent?

Ustinov(\cite{us}) already showed that 
they are not shift equivalent to each other
using ideas of combinatorial group theory.
We will give an alternate argument
using the canonical SFT covers.
We will compare the Bowen-Franks groups of
the SFT covers defined by period-2 orbits.
Since conjugacy preserves flow equivalence classes,
it suffices to show that 
there is no bijection between
the SFT covers defined by period 2 orbits
in $(X,g_1)$ and $(X,g_2)$, respectively,
such that the bijection respects 
the Bowen-Franks groups.

The presentations $(X,g_1)$ and $(X,g_2)$
are given by the following diagrams. 
\begin{figure}[ht]
\centering
\includegraphics[width=3cm]{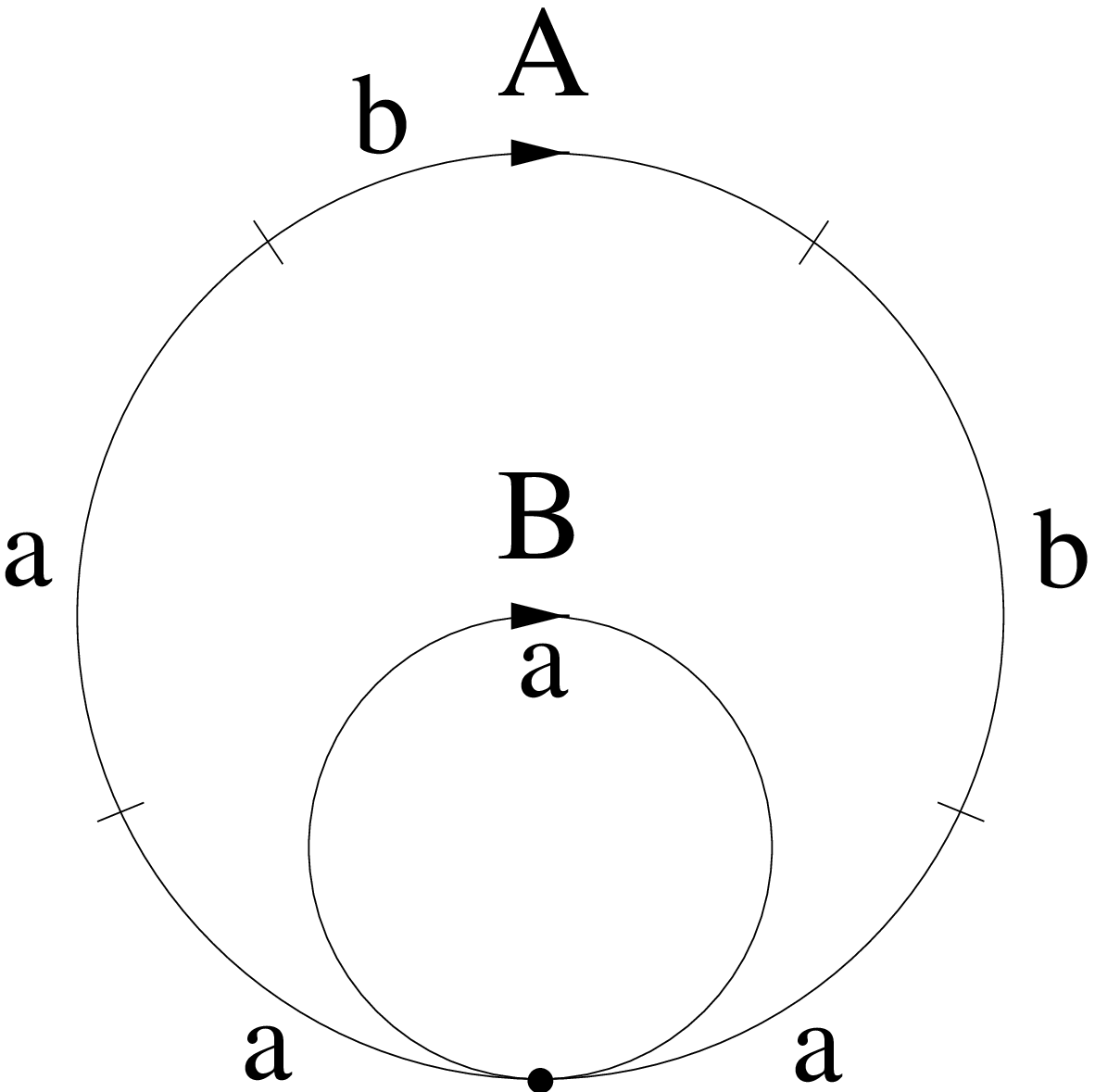}
\hspace{1.5cm}
\includegraphics[width=3cm]{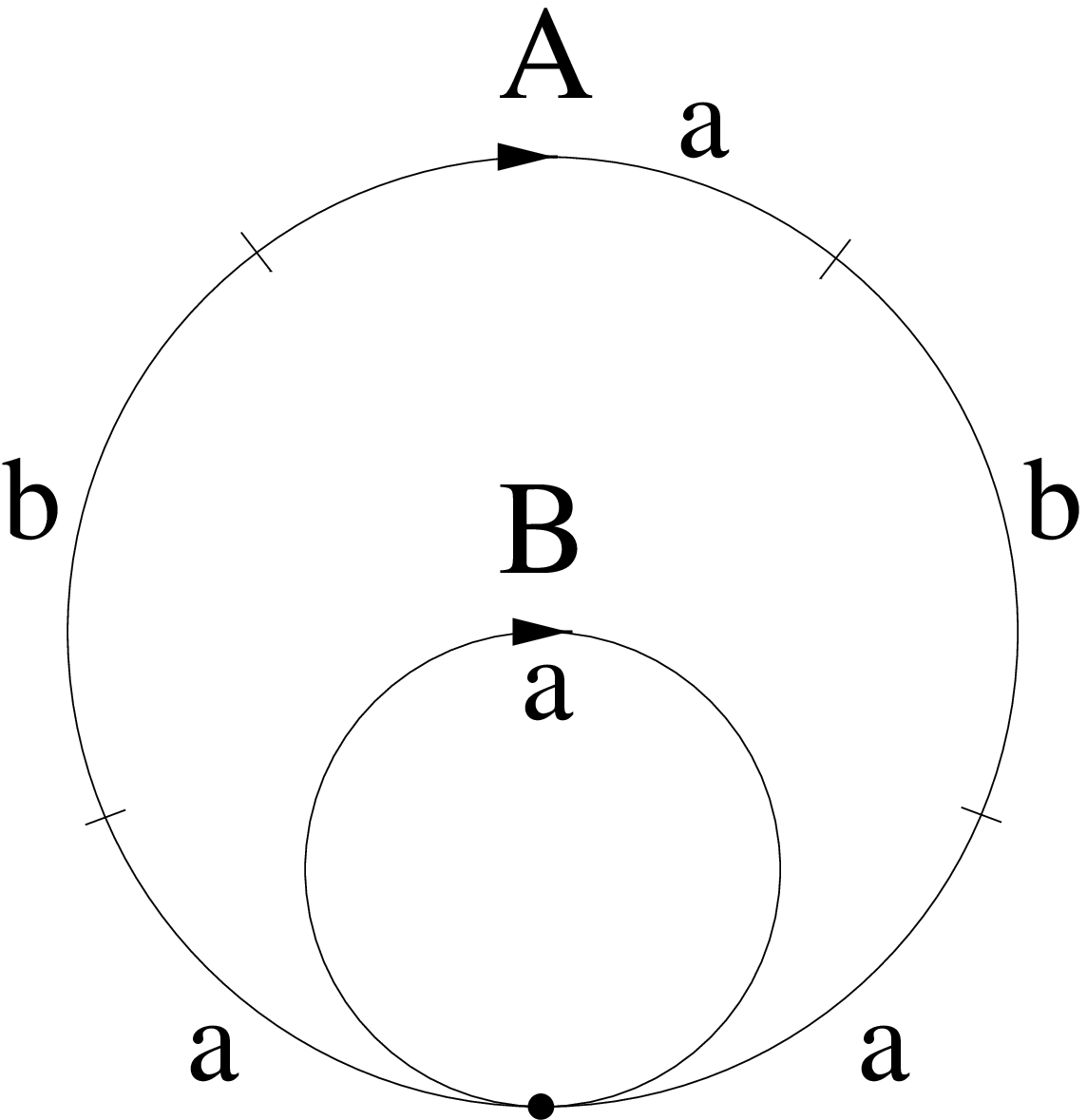}
\end{figure}
So $A$ has three subpaths $a_1,a_2,a_3$
which map to $A$ by $g_i$,
and two subpaths $b_1,b_2$ which map to $B$.

Each presentation has 
five period 2 orbits (excluding fixed points).
Let's denote $(a_i,a_j)$ as
the period $2$ orbit contained in
$a_i\cup a_j$
and $(b_k,B)$ as  the period $2$ orbit in
$b_k\cup B$.
For $(X,g_1)$,
the SFT covers defined by period 2 orbits are 
represented by the following matrices
\[
M(a_1,a_2)=
\begin{pmatrix} 0&1&0\\ 3&1&1\\ 5&3&1
\end{pmatrix}
,\,\,
M(a_1,a_3)=
\begin{pmatrix} 1&2&0\\ 1&1&1\\ 3&4&0 \end{pmatrix}
,\,\,
M(a_2,a_3)=
\begin{pmatrix} 2&3&0\\ 0&0&1\\ 4&5&0\end{pmatrix}
\]
\[
M(b_1,B)=
\begin{pmatrix}
2&1&1&1\\ 1&0&1&1\\ 1&0&0&0\\ 2&1&0&0
\end{pmatrix}
,\text{ and }\quad
M(b_2,B)=
\begin{pmatrix}
2&1&1&1\\  1&0&0&1\\ 1&0&0&0\\ 2&1&1&0
\end{pmatrix}
\]
and for $(X,g_2)$,
\[
N(a_1,a_2)=
\begin{pmatrix} 0&0&1\\ 3&1&1\\ 4&2&1\end{pmatrix}
,\,\,
N(a_1,a_3)=
\begin{pmatrix} 1&0&2\\ 2&1&0\\ 3&2&0 \end{pmatrix}
,\,\,
N(a_2,a_3)=
\begin{pmatrix} 1&3&1\\ 0&0&1\\ 2&4&1 \end{pmatrix}
\]
\[
N(b_1,B)=
\begin{pmatrix} 1&2&2\\ 1&1&2\\ 1&1&0 \end{pmatrix}
,\text{ and }\quad
N(b_2,B)=
\begin{pmatrix}
1&2&2\\ 1&0&1\\ 1&2&1 \end{pmatrix}
\]

We indicate the computation of $M(a_1,a_2)$
as an example.
The points  in the periodic orbit $(a_1,a_2)$
and the new graph $X_{(a_1,a_2)}$
are given in the following diagrams.
\begin{figure}[ht]
\centering
\includegraphics[width=3cm]{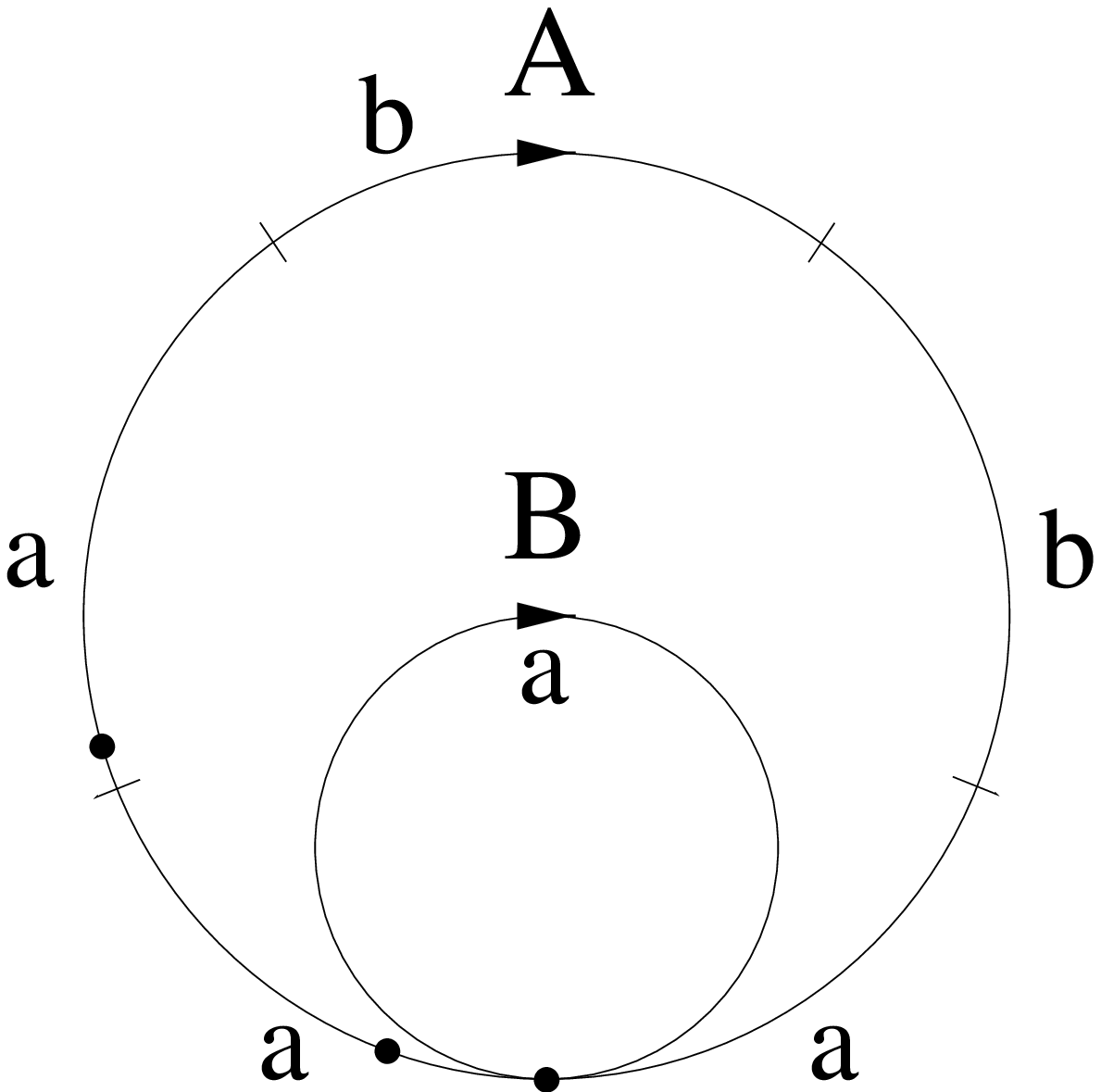}
\hspace{1.5cm}
\includegraphics[width=2.2cm]{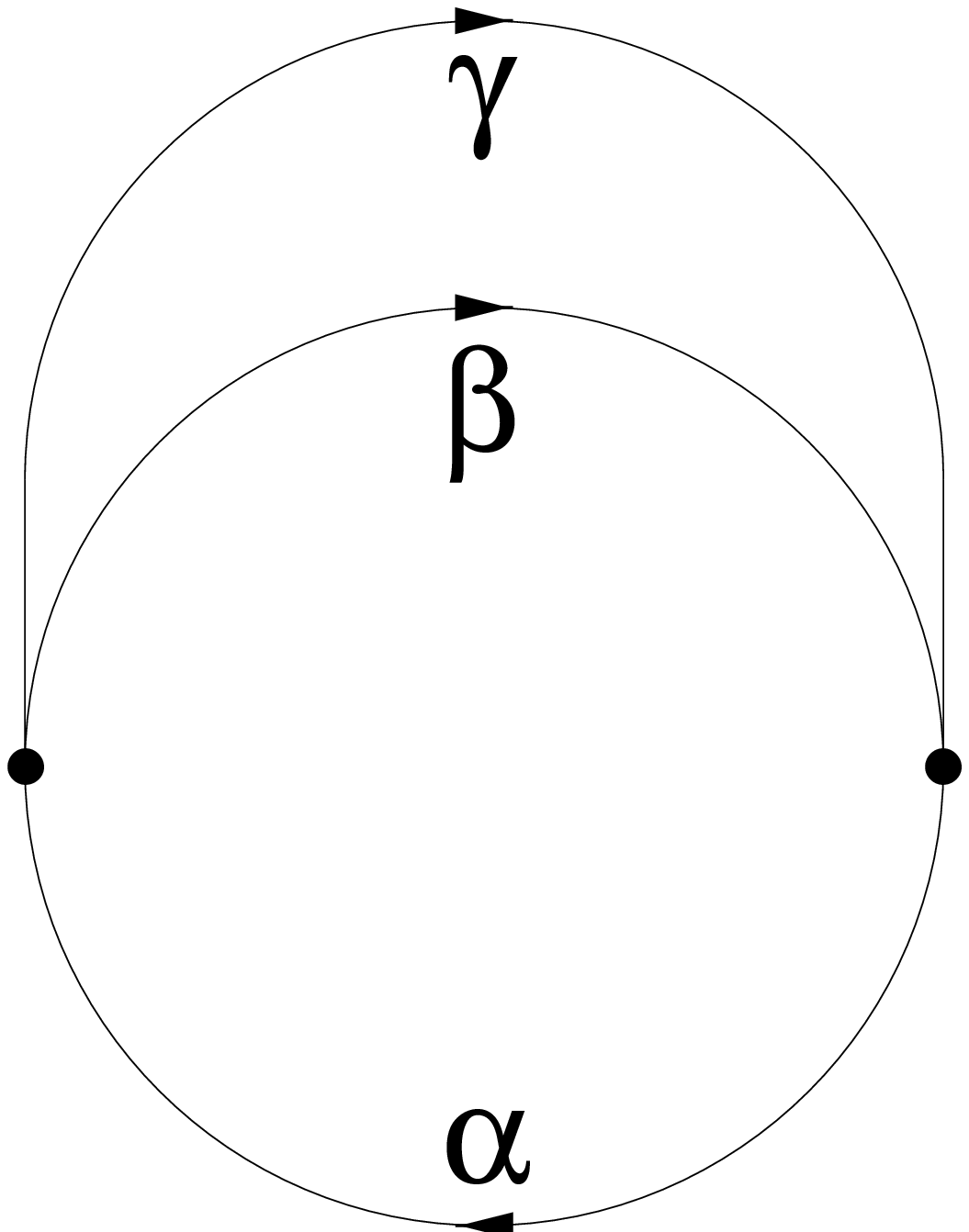}
\end{figure}
Let $p_1$ denote the path
from the branch point to the point in $a_1$,
$p_2$ the path from the point in $a_1$
to the point in $a_2$,
and $p_3$ the path from the point in $a_2$
to the branch point.
Let $\alpha=p_2$, $\beta=p_3p_1$,
and $\gamma=p_3bbp_1$.
Then the substitution rule
$\tilde{g_1}_{(a_1,a_2)}\colon
\mathcal{E}_{(a_1,a_2)}\to\mathcal{E}_{(a_1,a_2)}^*$
is given by
\[
\alpha\mapsto \beta,\quad
\beta\mapsto \alpha\gamma\alpha\beta\alpha,\quad
\gamma\mapsto \alpha\gamma\alpha\beta\alpha
                  \beta\alpha\beta\alpha.
\]

Use the Smith form to compute 
Bowen-Franks groups  (\cite[\S7.4]{lm}) 
of SFT covers defined by period 2 orbits.
Then it is not difficult to obtain
that in $(X,g_1)$, 
$M(a_1,a_3)$ and $M(b_1,B)$ have
$\mathbb{Z}_2\oplus \mathbb{Z}_4$, and
$M(a_1,a_2), M(a_2,a_3)$ and $M(b_2,B)$ have
$\mathbb{Z}_8$ 
as their Bowen-Franks groups.
And in $(X,g_2)$,
$N(a_1,a_2), N(a_1,a_3)$ and $N(a_2,a_3)$ have 
$\mathbb{Z}_2\oplus \mathbb{Z}_4$, 
and $N(b_1,B), N(b_2,B)$ have $\mathbb{Z}_8$.
So the number of SFT covers which have the same
Bowen-Franks groups are different, 
and $\overline{g}_1$ is not conjugate to
$\overline{g}_2$.
Therefore  $(X,g_1)$ is not shift equivalent to
$(X,g_2)$.
\end{example}

\appendix
\section{One-sided SFT}
Suppose that  $(X,f)$ and $(Y,g)$ are
presentations of solenoids 
which are shift equivalent to each other
by a shift equivalence $r\colon X\to Y$ and
$s\colon Y\to X$.
Assume that $\mathcal{O}$ is a finite invariant subset
of $(X,f)$ and $\mathcal{O}^\prime=r(\mathcal{O})$,
and denote $(\Sigma_{\mathcal{O}},\sigma_{\mathcal{O}})$
and $(\Sigma_{\mathcal{O}^\prime},
              \sigma_{\mathcal{O}^\prime})$
as the SFT  covers of $(\overline{X},\overline{f})$
and $(\overline{Y},\overline{g})$ defined by
$\mathcal{O}$ and $\mathcal{O}^\prime$ respectively.
Then $(\Sigma_{\mathcal{O}},\sigma_{\mathcal{O}})$ and
$(\Sigma_{\mathcal{O}^\prime},
              \sigma_{\mathcal{O}^\prime})$
are defined by 
nonnegative integer matrices $M_{\mathcal{O}}$ and
$M_{\mathcal{O}^\prime}$, respectively.
And we can make {\bf one-sided} SFTs
$(S_{\mathcal{O}},\sigma_{\mathcal{O}})$ and
$(S_{\mathcal{O}^\prime},\sigma_{\mathcal{O}^\prime})$  
from $M_{\mathcal{O}}$ and $M_{\mathcal{O}^\prime}$,
respectively.

The purpose of this appendix is to give an example
in which the one-sided SFTs
$(S_{\mathcal{O}},\sigma_{\mathcal{O}})$ and
$(S_{\mathcal{O}^\prime},\sigma_{\mathcal{O}^\prime})$ 
are not conjugate.

\begin{example} \label{a.2}
Let $(Y_{\{p\}},g_{\{p\}})$
be as in Example \ref{4.4}.
Let $Y_1$ be a graph such that
$\mathcal{E}_{Y_1}=\{1,2,3,4,5,6,7,8\}$
with a graph map $g_1\colon Y_1\to Y_1$
defined by 
\[
1\mapsto 1234357,\text{ }\text{ }
2,4,5,6\mapsto 8,\text{ }\text{ }
3\mapsto 167, \text{ }\text{ }
7\mapsto 1,
\text{ and }
8\mapsto 678.
\]
The presentation $(Y_{\{p\}},g_{\{p\}})$ and
the graph $Y_1$ are given in the following diagrams.
In $Y_1$, the fixed point labeled $p$ is
the terminal point of the edge 8
and the initial point of the edge 1.
\begin{figure}[!h]
\centering
\includegraphics[height=3.3cm]{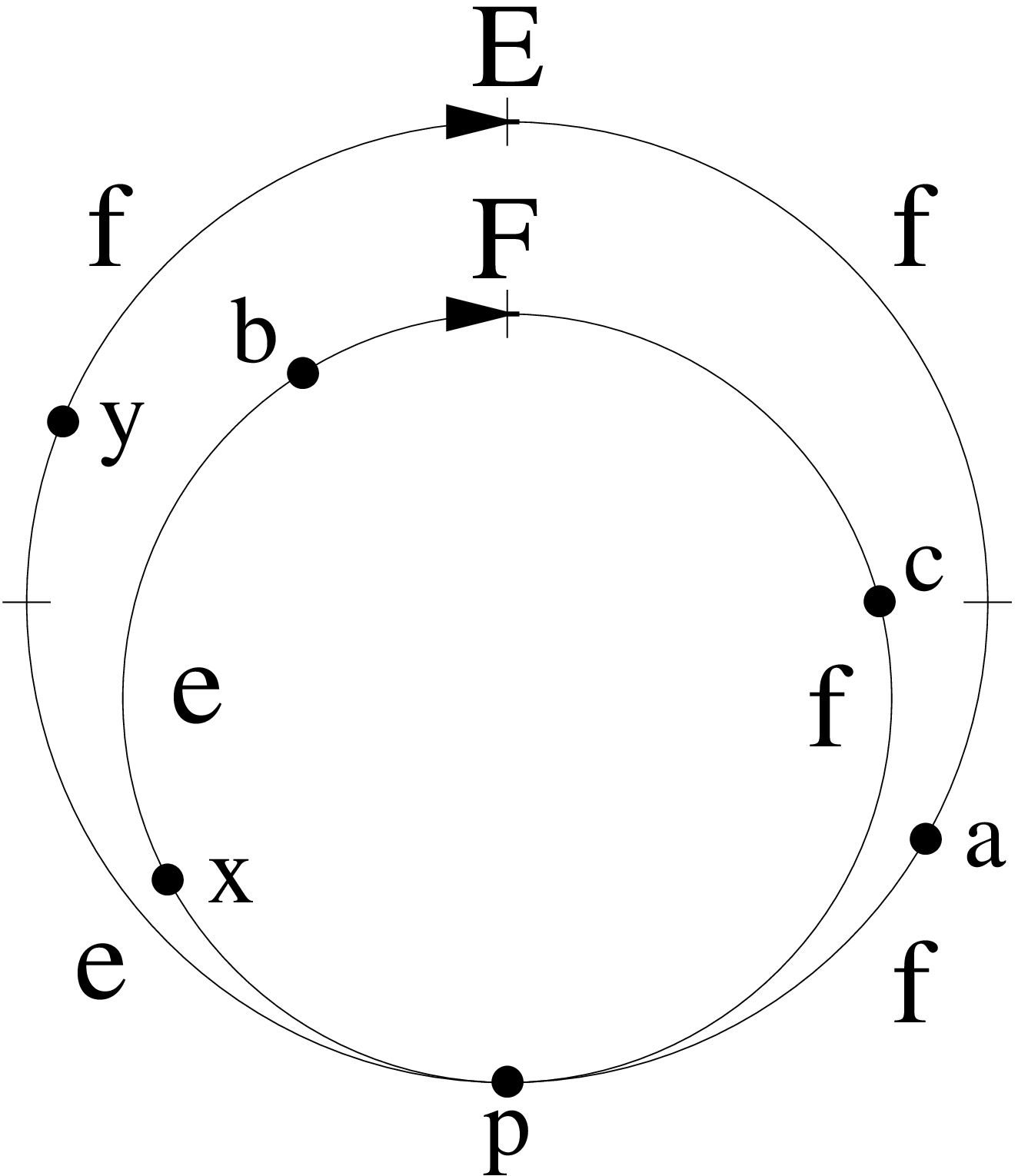}
\hspace{1cm}
\includegraphics[height=2.8cm]{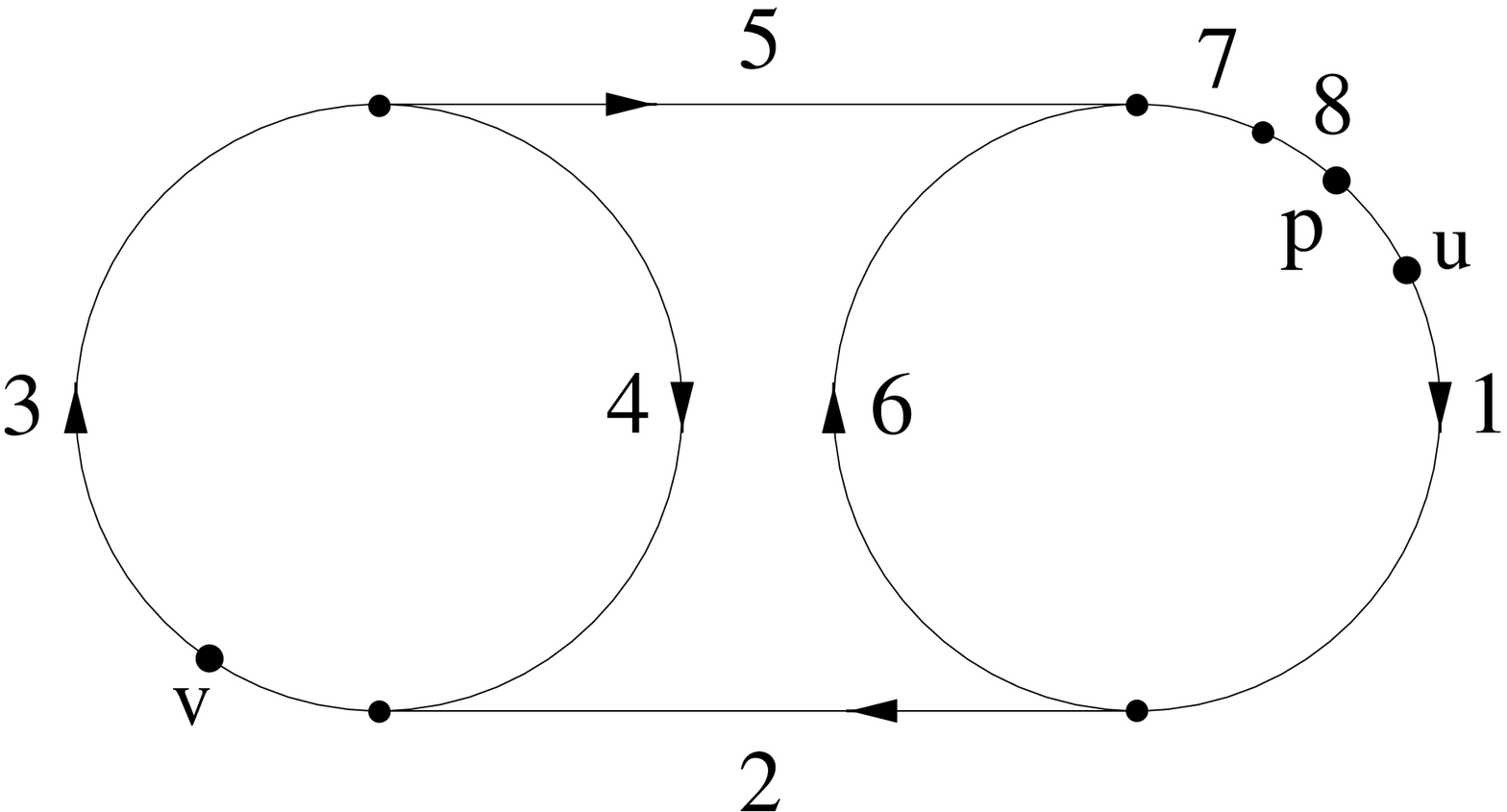}
\end{figure}

We will define a lag 1 shift equivalence  
of $g_{\{p\}}$ and $g_1$ by graph maps $r,s$
under which the points labeled $p$ correspond.
The points labeled $x,y$ in $Y_{\{p\}}$
and $u,v$ in $Y_1$ form period two orbits of
$g_{\{p\}}$ and $g_1$, respectively. 
The map $r$ will send the points $x,y$ to
the points $u,v$. 
But, we will see that the covering SFTs associated 
to $\{x,y\}$ and $\{u,v\}$ are not conjugate 
as one-sided shifts.

The points $a,b$ in $Y_{\{p\}}$ is a period two orbit,
and the point $c$ is the unique  point in the edge $F$
such that $g_{\{p\}}(c)=b$.
Let $e_1$ be the path from $p$ to $a$,
$e_2$ the path from $a$ to $p$,
$f_1$ the path from $p$ to $b$,
$f_2$ the path from $b$ to $c$,
and $f_3$ the path from $c$ to $p$.
Define $r\colon Y_{\{p\}}\to Y_1$ and
$s\colon Y_1\to Y_{\{p\}}$ by
\begin{align*} 
\tilde{r} &\colon e_1\mapsto 1234357,\  e_2\mapsto 8, \ 
f_1\mapsto 1,\  f_2\mapsto 67, \  f_3\mapsto 8 \\  
\tilde{s}&\colon 1\mapsto e_1, \,\  2,6\mapsto e_2,
\, 3\mapsto f_1f_2,\  \, 4,5\mapsto f_3,  
\,\  7\mapsto f_1, \,\  8\mapsto f_2f_3.
\end{align*}
Then $ \tilde{r}$ is given by
$e\mapsto 12343578$ and $f\mapsto 1678$, and 
\begin{align*}
\tilde{s}\circ \tilde{r}&\colon
e\mapsto
e_1e_2f_1f_2f_3f_1f_2f_3f_1f_2f_3=efff,
\text{ and }
f\mapsto
e_1e_2f_1f_2f_3=ef\\
\tilde{r}\circ \tilde{s}&\colon
1\mapsto 1234357,\ 
2,4,5,6\mapsto 8, \ 
3\mapsto 167, \ 
7\mapsto 1, \ 
\text{ and }
8\mapsto 678.
\end{align*}
Therefore we have $s\circ r=g_{\{p\}}$ and $r\circ s=g_1$,
and $(Y_{\{p\}},g_{\{p\}})$ is 
shift equivalent to $(Y_{1},g_{1})$ by $r$ and $s$.

Now $\mathcal{P}_{\{x,y\}}=\{\alpha, \beta, \gamma\}$ 
where $\alpha$ is the path from $y$ through $p$ to $x$,
$\beta$ is the circle $F$ based at $x$, and
$\gamma$ is the path from $x$ through $p$ to $y$.
And $\mathcal{P}_{\{u,v\}}
=\{\delta, \epsilon, \zeta, \eta\}$ where
$\delta$ is the path from $u$ through the path $2$ to $v$,
$\epsilon$ is the circle $34$ based at $v$,
$\zeta$ is the path from $v$ through the path $5$ to $u$,
and $\eta$ is the circle $1678$ based at $u$.
The wrapping rules
$\tilde{g}_{\{x,y\}}\colon
\mathcal{P}_{\{x,y\}} \to \mathcal{P}_{\{x,y\}}^*$
and $\tilde{g}_{\{u,v\}}\colon
\mathcal{P}_{\{u,v\}} \to \mathcal{P}_{\{u,v\}}^*$
are given by
\begin{align*}
\tilde{g}_{\{x,y\}}&\colon
\alpha\mapsto \beta\beta\gamma,\quad
\beta\mapsto \alpha\gamma,\quad
\gamma\mapsto \alpha\gamma\alpha, \text{ and}\\
\tilde{g}_{\{u,v\}}&\colon
\delta \mapsto \epsilon \zeta,\quad
\epsilon \mapsto \eta,\quad
\zeta \mapsto \eta \eta \delta,\quad
\eta \mapsto \epsilon \zeta \eta \delta.
\end{align*}
Therefore
the SFT covers associated to 
$\{x,y\}$ and $\{u,v\}$ are
presented by the following matrices.
\[
M_{\{x,y\}}=
\begin{pmatrix}
0&2&1\\
1&0&1\\
2&0&1
\end{pmatrix}
,\qquad
M_{\{u,v\}}=
\begin{pmatrix}
0&1&1&0\\
0&0&0&1\\
1&0&0&2\\
1&1&1&1
\end{pmatrix}
\]

For 
\[
R=
\begin{pmatrix}
0&1&1&0\\
0&0&0&1\\
1&0&0&1
\end{pmatrix}
\text{ and }
S=
\begin{pmatrix}
1&0&0\\
0&1&0\\
0&1&1\\
1&0&1
\end{pmatrix}
, 
\]
we have 
$M_{\{x,y\}}=RS $ and $M_{\{u,v\}}=SR $, and 
the two-sided SFT covers are 
topologically conjugate (as guaranteed by 
Theorem \ref{5.aa}). 
On the other hand, 
Williams showed \cite{w4} that  one-sided 
SFTs are conjugate if and only if they have 
the same total column amalgamation matrix; and
we   see that 
$M_{\{x,y\}}$ is its own total column amalgamation, 
whereas the 
total column amalgamation of $M_{\{u,v\}}$ is
\[
\begin{pmatrix}
0&1&0\\
1&0&3\\
1&1&1
\end{pmatrix} .
\]
Therefore these one-sided SFTs are not conjugate.
\end{example}


\section{Relation between Williams' definition and 
topological definition of $1$-solenoids}

Let $(X,f)$ be a presentation of
a solenoid in our topological sense, that is, 
assume that it satisfies the Axioms 0-5 of Section 2.
In this appendix,  
we will give $X$ the differentiable structure 
of a branched $1$-manifold,
with respect to which $f$ will be an immersion.
This will show that 
our topological $1$-solenoids are exactly the 
systems obtained from Williams' differentiable 
$1$-solenoids by forgetting the differentiable 
structure. 

The precise definition of
an $n$-dimensional branched manifold (\cite{w3})
is necessarily somewhat complicated and subtle.
In our one-dimensional situation,
we will attach to $(X,f)$
a more simple and special structure,
from which an interested reader can easily derive
the full immersion of a branched manifold
apparatus laid out in \cite{w3}.
There is some  overlap in our ideas and
those used by Williams for his Realization Theorem 
\cite[7.6]{w2}.

We will cover the graph $X$ with open sets 
$V_1, \dots , V_k$.
Each $V_i$ will be an interval or
a union of intervals intersecting at a vertex.
There will be attached maps ({\it charts})
$\pi_i\colon V_i\to I_i$,
where the $I_i$ are open intervals.
When $V_i$ is an interval, 
$\pi_i$ will be a homeomorphism;
in the other case
$V_i$ will be a union of intervals the restriction of 
$\pi_i$ to each of which is a homeomorphism.
Let $V_{ij}= V_i\cap V_j$.
Whenever $i\neq j$ and $V_{ij}$ is nonempty, 
the restriction of the map $\pi_i$ to $V_{ij}$
will be invertible and
the `{change-of-coordinates}' map
\[
\pi_{i,j} \colon  \pi_i(V_{ij})\to \pi_j(V_{ij})
\text{ given by } x\mapsto \pi_j\circ \pi_i^{-1}x
\]
will be invertible and affine. 
Finally, for each pair $i,j$ such that
$f(V_i)\cap V_j$ is nonempty,
there will be an invertible  affine map 
$f_{i,j}\colon I_i \to I_j$ such that 
the restriction of $f$ to $V_i\cap f^{-1}(V_j)$ is 
a lift of $f_{i,j}$, that is to say,  
$f_{i,j}\circ \pi (x) = \pi_j\circ f(x)$ 
for all $x$ in $V_i\cap f^{-1}(V_j)$. 
 
To define the open sets $V_i$,
fix $N$ such that  
each point of $X$ has a neighborhood whose image 
under $f^N$ is an interval such that 
the interior of every interval with endpoints     
in  $f^{-1}(V)$ contains a point in 
$f^{-N}(V)\setminus f^{-1}(V)$. 
The sets $V_i$ will be of two sorts. 
First, the complement of $f^{-N}(V)$ in $X$
will be a disjoint union of open intervals, 
and each of these will be one of the sets $V_i$.
Second, at each point $y$ in $f^{-N}(V)$,
we pick a connected open neighborhood $U_y$
such that the $U_y$ are disjoint;
after some shrinking,
these will be the remaining $V_i$. 
The complement in $U_y$ of $y$ will be the union 
of a collection $\mathcal J(y)$ of
disjoint open intervals $J_{y,t}$
(if $y$ is not a vertex, then there are just two);
we pass to a shrunken set of  $V_i$ for which  
the images of any pair $J_{y,t},J_{y',t'}$ 
are equal or disjoint.
This completes the description of the $V_i$. 
  
To describe the charts,
we will use the following result of Williams.
Williams proved it for elementary presentations,
but his proof works for the non-elementary 
case and it works with our topological axioms.  

\begin{lemma}[{\cite[6.2]{w2}}]\label{b.3}
There exist a unique measure $\mu$ on $X_1$ and
a unique real number $\lambda>1$ such that
$\mu (X_1)=1$ and
$\mu \bigl(f_1(I)\bigr)=\lambda I$
for every small interval contained in an edge.
\end{lemma}

First we use this measure $\mu$ to define 
the charts $\pi_i$ when $V_i$ is an interval:  
identifying $V_i$ with 
$(0,1)$, we define $\pi_i(x) = \mu (0,x)$.  
Next, for each $y$ in $f^{-N}(V)$, 
define an equivalence relation on $\mathcal J(y)$ 
by declaring  two intervals to be equivalent
if their images under $f^N$ are equal. 
This divides each $\mathcal J(y)$ into two 
equivalence classes; making an arbitrary choice, 
designate one class as $\mathcal L(y)$ and 
the other as $\mathcal R(y)$.
 
Now we can describe the chart $\pi_i$ for $V_i=U_y$.   
Identify each $J_{y,t}$ with the interval $(0,1)$,
with $0$ corresponding to $v$, 
and  for $x$ in $J_{v,i}$ define 
\[
\pi_v(x)=\begin{cases}
\mu (0,x) \qquad &\text{if } J_{v,i}\in R(v)\\   
-\mu (0,x)       &\text{if } J_{v,i}\in L(v)\ .  
\end{cases}
\]
Finally, define $\pi_i(y)=0$.
This completes the definition of the charts $\pi_i$. 
It is easy to verify that
the change-of-coordinate maps $\pi_{i,j}$ are 
one-to-one and affine as claimed. 

It remains to see that
$f$ is locally the lift of affine maps as claimed.
Suppose that the set
$V_{i,j}=V_i\cap f^{-1}(V_j)$ is nonempty.  
If $V_i$ is not one of the $U_y$,
then by choice of $N$,
$f(V_i)\cap V_j$ is an interval,
and by choice of $N$ if $V_j=U_y$,
then this interval is entirely contained in
one of the intervals $J(y,t)$; 
so the required affine map $f_{i,j}$
(with multiplicative constant $\lambda$ or $-\lambda$)
exists. 
Similarly the  required $f_{i,j}$ exists if  
$V_i=U_y$ and $f(y)$ is not a vertex.   

Next, assume  $V_i=U_y$ and $f(y)$ is a vertex.
For each $J_{y,t}$,
let $J'_t$ denote the unique member of
$\mathcal J(f(y))$ intersecting $f(J_{y,t})$.  
First suppose that $y$ is not a vertex, 
so $\mathcal J(y)= \{J_{y,1},J_{y,2}\}$. 
By the Nonfolding Axiom, 
the restriction of $f^{N+1}$ to $V_i$ is 
locally injective at $y$, 
and therefore $J'_1$ and $J'_2$ must be in 
different $\mathcal L / \mathcal R$ classes.  
Therefore the required affine map 
$f_{i,j}$ exists,  with multiplicative 
constant $\lambda$ or $-\lambda$.

Finally, suppose  $V_i=U_y$,  
$y$ is a vertex and (therefore) $f(y)$ is a vertex.   
We claim that for any $J_1$ in $\mathcal L(y)$ 
and $J_2$ in $\mathcal R(y)$, the intervals 
$J'_1$ and $J'_2$ cannot be in the same equivalence 
class $\mathcal L(f(y))$ or $\mathcal R(f(y))$. 
Given the claim, we can define the required 
map $f_{i,j}$ as $x\mapsto \lambda x$
(if the $\mathcal L, \mathcal R$ 
equivalence class labeling is respected) 
or $x\mapsto -\lambda x$
(if the labeling is reversed).
So it remains to prove the claim. 

Suppose the claim is false - suppose there are   
$J_1\in \mathcal L(y)$ and $J_2 \in \mathcal R(y)$ 
with $J'_1$ and $J'_2$ both in $\mathcal L(f(y))$
or both in $\mathcal R(f(y))$.
Because $f^N(J'_1)=f^N(J'_2)$ 
it follows that $f^{N+1}$ is not locally injective 
at $y$ on the interval which is the union of 
$\{y\}$, $J_1$ and $J_2$.
Because the image of under $f^N$
of an interval $J_{y,t}$ depends only 
on the  class $\mathcal L(y)$ or $\mathcal R(y)$ 
to which it belongs, it follows that 
the restriction of $f^{N+1}$ to  any open 
interval containing $y$ is not locally one-to-one. 
Now pick a positive integer $k$ and a point $w$ 
which is not a vertex such that $f^k(w)=y$. 
Then $f^{k+N+1}$ is not locally injective at $w$. 
This contradicts the Nonfolding Axiom, and 
concludes the proof. We record the result as 
the following Proposition.  

\begin{proposition} 
Suppose $(X,f)$ is a presentation of a 
topologically defined $1$-solenoid.
Then $X$ can be given a differentiable structure 
making it a branched manifold on which $f$ 
is an immersion presenting a Williams solenoid.  
\end{proposition}

\begin{example}\label{b.e}
Let $X$ be given by the following graph,
\begin{figure}[!h]
\centerline{\scalebox{.32}{\includegraphics{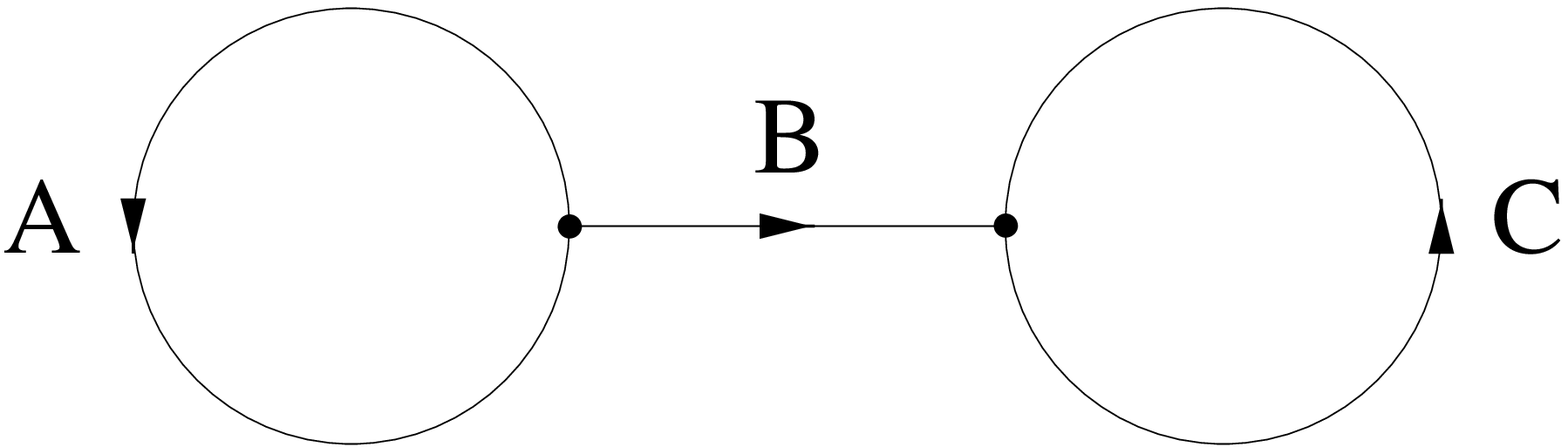}}}
\end{figure}
and $f\colon X\to X$
is given by
$a\mapsto b^{-1}ab,\text{ }b\mapsto cb^{-1}a,
\text{ }c\mapsto bcb^{-1}$.

If we redraw $X$ as the following graph,
then $X$ is a branched manifold,
and $(X,f)$ is a presentation
in the sense of Williams.
\begin{figure}[!h]
\centerline{\scalebox{.29}{\includegraphics{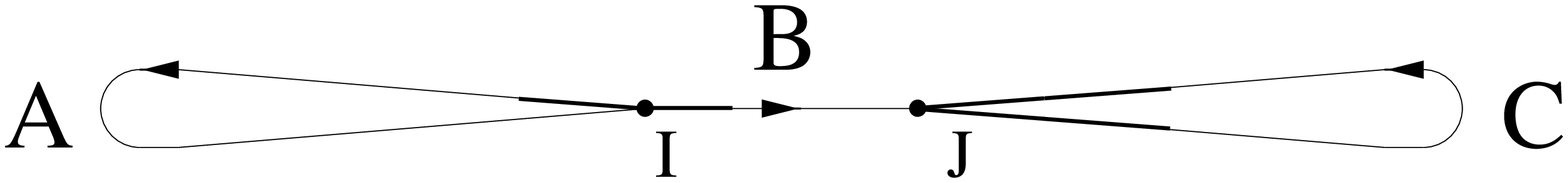}}}
\end{figure}
\end{example}

In contrast to the solenoid case,
the following example shows that
the set of Williams' branched solenoids
is a proper subset of the topologically defined 
branched solenoids.

\begin{example}\label{b.ee}
Let $X$ be as in the previous example,
and $g\colon X\to X$ is given by
\[
a\mapsto cb^{-1}ab,\text{ }b\mapsto c^{-1}b^{-1},
\text{ }c\mapsto a.
\]
Then $(X,g)$ satisfies all Axioms
except for the Flattening Axiom.
So $(X,g)$ is a presentation of a 
branched solenoid
according to the topological definition.  

At each branch point of $X$,
there are three choices of differentiable structure,
that is,
three arcs are parallel to each other
or two arcs are parallel and the other
arc is 180 degrees to these two arcs.
The second graph in Example \ref{b.e}
and Figure \ref{xcv21} show three possible
differentiable structures at the left branch point
when the circle $c$ is fixed.
And similar differentiable structures
can be given to the right branch point.

In each choice of differential structures
at both branched points,
it is not difficult to find
a smooth interval which is mapped to
a non-smooth interval homeomorphically
by $g\colon X\to X$.
For example, if $X$ is the second graph
in Example \ref{b.e},
then the interval $I$ is mapped to
the interval $J$ by $g$ homeomorphically.
Hence $(X,g)$ cannot be a presentation
of a branched solenoid according to
Williams' definition.
\begin{figure}[ht]\label{xcv21}
\centering
\begin{minipage}[c]{0.425\textwidth}
\centering \includegraphics[height=1.5cm]{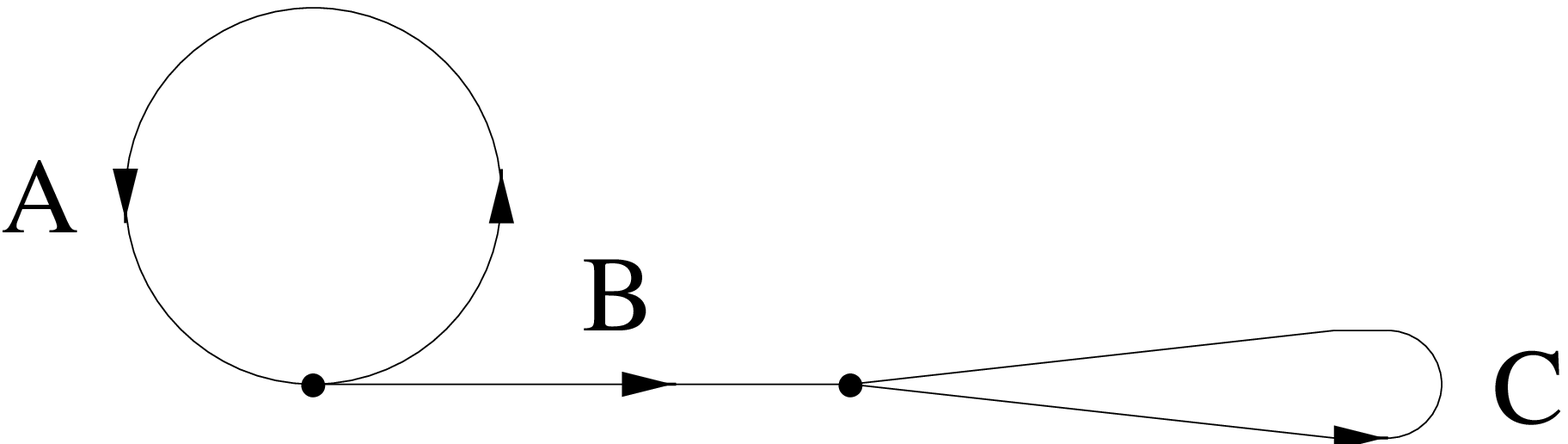}
\end{minipage}
\hskip 0.6cm
\begin{minipage}[c]{0.425\textwidth}
\centering \includegraphics[height=1.5cm]{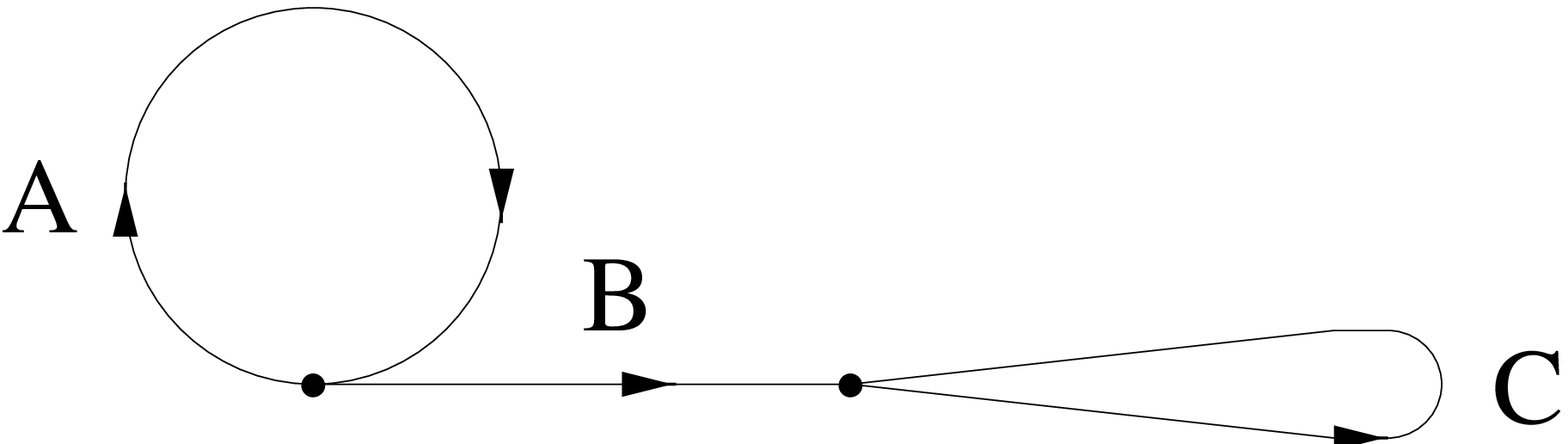}
\end{minipage}
\caption{Other differentiable structures at
the branch point of $A$ when the circle $c$ is fixed.}
\end{figure}
\end{example}

\begin{ack}
This paper will be  part of my Ph.D. thesis at UMCP.
I would like to express deep gratitude
to my advisor Dr.{\;}Mike Boyle
for his encouragement and many useful discussions.
\end{ack}

\bibliographystyle{amsplain}

\end{document}